\documentclass{article}
\usepackage{relsize,calc,bigstrut,afterpage,rotfloat,multirow,array,amsmath,amsthm,amssymb,mathrsfs}
\usepackage{hyperref}
\usepackage[plain]{algorithm}
\usepackage{algpseudocode}
\usepackage[all]{xy}

\newcommand{\Star}{\ast}

\newtheorem{theorem}{Theorem}[section]
\newtheorem{definition}[theorem]{Definition}
\newtheorem{proposition}[theorem]{Proposition}
\newtheorem{corollary}[theorem]{Corollary}
\newtheorem{lemma}[theorem]{Lemma}
\newtheorem{fact}[theorem]{Fact}

\newtheorem*{remark}{Remark}
\newtheorem*{problem}{Open Problem}
\newtheorem*{question}{Question}

\newcommand{\minisection}[1]{\vspace{0.15cm}\noindent {\bfseries #1.}}
\newenvironment{example}{\begin{proof}[\textbf{\upshape Example}]}{\end{proof}}

\newcommand{\A}{\mathscr{A}}
\newcommand{\B}{\mathscr{B}}

\newcommand{\Q}{\mathcal{Q}}
\renewcommand{\P}{\mathcal{P}}
\newcommand{\SN}{\mathscr{N}}
\newcommand{\SP}{\mathscr{P}}
\newcommand{\QP}{(\Q,\P)}

\newcommand{\Pow}{\textit{Pow}}
\newcommand{\cl}{\mathrm{cl}}

\newcommand{\opts}{\mathrm{opts}}
\newcommand{\<}{\langle}
\renewcommand{\>}{\rangle}

\newcommand{\PPos}{$\SP$-position}
\newcommand{\NPos}{$\SN$-position}

\newcommand{\CGSuite}{{\itshape\textsmaller{C\kern- 0.06emG\kern-0.1emS\kern-0.04em}uite}}
\newcommand{\sh}{\textrm{\protect\raisebox{1pt}{\tiny \#}}}

\font\Chess=chess15
\input macros.tex

\begin{document}

\bibliographystyle{abbrv}

\title{Mis\`ere Quotients for Impartial Games}

\author{Thane E. Plambeck\\2341 Tasso St.\\Palo Alto, CA 94301 \and Aaron N. Siegel\\Institute for Advanced Study\\1 Einstein Drive\\Princeton, NJ 08540}



\date{\today}
\maketitle

\vspace{-0.5cm}

\begin{center}
\itshape \small Dedicated to Richard K. Guy on the occasion of his $\textit{90}^\textit{th}$ birthday, 30 September 2006.
\end{center}

\vspace{0.15cm}

\begin{abstract}
We announce mis\`ere-play solutions to several previously-unsolved combinatorial games.  The solutions are described in terms of \emph{mis\`ere quotients}---com\-mu\-ta\-tive monoids that encode the additive structure of specific mis\`ere-play games.  We also introduce several advances in the structure theory of mis\`ere quotients, including a connection between the combinatorial structure of normal and mis\`ere play.
\end{abstract}

\section{Introduction}

\label{section:introduction}

\begin{figure}
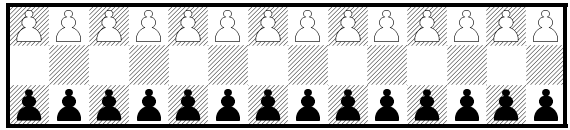

\centering

\makebox{
\vbox{\offinterlineskip%
   \hrule height1pt%
   \hbox{\vrule width1pt\Chess%
         \vbox{\hbox{OPOPOPOPOPOPOP}%
               \hbox{0Z0Z0Z0Z0Z0Z0Z}%
               \hbox{opopopopopopop}}%
         \vrule width1pt}%
   \hrule height1pt}%
 }
\caption{The starting position in $3 \times 14$ Dawson's Chess.}
\label{figure:dawsonschess}
\end{figure}

In 1935, T.~R.~Dawson, the prolific composer of fairy chess problems, invented a little game now known as \emph{Dawson's Chess} \cite{dawson_1935}.  The game is played on a $3 \times n$ chessboard, with equal facing rows of black and white pawns.  The $3 \times 14$ starting position is shown in Figure~\ref{figure:dawsonschess}.  The pieces move and capture as ordinary chess pawns, and captures, if available, are mandatory.  The mandatory-capture rule ensures that each pawn will advance at most
one rank before being captured or blocked by an opposing pawn.  
Thus the players will eventually run out of moves, and whoever makes the last move \emph{loses}.

Dawson invested considerable effort in finding a perfect strategy (on arbitrarily large boards), but he was hampered by his unfortunate choice of winning condition.  Had he asserted, instead, that whoever makes the last move \emph{wins}, he might have solved the game. 
Instead, the problem has remained open for over~70 years.

This is no accident: combinatorial games with \emph{mis\`ere play}---in which the last player to move loses---tend to be vastly more complicated than their normal-play (last-player-winning) counterparts.  Many authors have attempted to solve Dawson's Chess, and a variety of related mis\`ere-play games, only to be met with fierce resistance.

In this paper, we present strategies for several previously-unsolved mis\`ere-play games.  The solutions were obtained using a combination of new theoretical advances and powerful computational techniques.  These advances also shed light on the question of just \emph{why} mis\`ere play is so much harder than normal play.  Although Dawson's Chess remains open, we have pushed the analysis further than ever before; and there is significant hope that our techniques can be improved to obtain a full solution.

\subsection*{Impartial Combinatorial Games}

A combinatorial game is \emph{impartial} if both players have exactly the same moves available at all times.\footnote{Although Dawson's Chess appears to distinguish between Black and White, the reader might wish to verify that this is illusory: the mandatory-capture rule guarantees that their moves are, in fact, equivalent.}  A complete theory of \emph{normal-play} impartial games was convincingly developed in the 1930s, independently by R.~P.~Sprague and P.~M.~Grundy \cite{sprague_1935, sprague_1937}.  They showed that every position in every impartial game is equivalent (in normal play) to a single-heap position in the game Nim.

In the 1956 issue of the \emph{Proceedings of the Cambridge Philosophical Society}, there appeared two papers, back-to-back, that highlighted the sharp differences between normal and mis\`ere play.  The first was a seminal article by Richard K.~Guy and Cedric A.~B.~Smith, \emph{The $G$-Values of Various Games} \cite{guy_1956}, in which the authors introduced a broad class of impartial games known as \emph{octal games}.  Guy and Smith applied the Sprague--Grundy theory to obtain complete normal-play solutions to dozens of such games, including Dawson's Chess.\footnote{In fact, in 1947 Guy succeeded Dawson as endgames editor of the \emph{British Chess Magazine}, and he credits Dawson's Chess as a motivating influence: ``It may be of historical interest to note that Dawson showed the problem to me around 1947.  Fortunately, I forgot that Dawson proposed it as a losing game, was able to analyze the normal-play version and rediscover the Sprague--Grundy theory.'' \cite{guy_2006pc}}  Their analysis was such an enormous success that it left little to be discovered about octal games in normal play.

The second article was by Grundy and Smith \cite{grundy_1956}, and concerned impartial games with mis\`ere play.  In contrast to the Guy--Smith paper, it has a somewhat dispiriting introduction:

\begin{quotation}
\small
Various authors have discussed the `disjunctive compound' 
of games with the last player winning (Grundy~\kern-.6pt\cite{grundy_1939}; 
Guy and Smith~\kern-.6pt\cite{guy_1956}; Smith~\kern-.6pt\cite{smith_unpub}).  We attempt here to 
analyse the disjunctive compound with the last player losing, 
though unfortunately with less complete success \dots.
\end{quotation}

At the source of the ``less complete success'' was the apparent lack of a natural analogue of the
Sprague--Grundy theory in mis\`ere play---in particular, the almost negligible canonical simplification of most mis\`ere-play
impartial game trees.  Grundy and Smith speculated that improvements might be lurking, but their efforts to produce additional simplication rules were unsuccessful.  Finally, in 1976 Conway proved through an intricate argument \cite[Theorem 77]{conway_1976} that no such improvements exist.  Thirty years later, he offered the following remarks.

\begin{quotation}
\small
The result I am most proud of in mis\`ere impartial combinatorial games is that the Grundy and Smith reduction
rules are in fact the only ones available in general in the global semigroup of all mis\`ere games \dots.  Grundy quite specifically did not
believe that this would be the case.~\cite{conway_2006pc}
\end{quotation}

Conway's result shows that the complications observed by Grundy and Smith are intrinsic.  Yet for many mis\`ere games, including many octal games, there is a path forward.  The Grundy--Smith equivalence can be localized to the set of all positions that arise in the play of a particular game.  This associates to every game $\Gamma$ a certain commutative monoid~$\Q$, the \emph{mis\`ere quotient} of $\Gamma$.  Together with a small amount of additional information, one can recover from~$\Q$ a perfect winning strategy for $\Gamma$.  The monoid~$\Q$ therefore serves as a \emph{local} mis\`ere analogue of the Sprague--Grundy theory.

The main ideas behind the quotient construction are described in~\cite{plambeck_2005,plambeck_200X}, and we'll review them in the next section; but they only set the stage for the goals of this paper:
\begin{enumerate}
\item To record what's been discovered about the mis\`ere play of particular games, in the spirit of Guy and Smith; 
\item To develop a structure theory for mis\`ere quotients; and
\item To identify new problems and areas of interest in the structure of mis\`ere play.
\end{enumerate}

Henceforth, we will assume that the reader is familiar with the classical Sprague--Grundy theory.  To a lesser extent, we will also assume familiarity with the canonical theory of mis\`ere games.  See \emph{Winning Ways}~\cite[Chapters~4~and~13]{berlekamp_1982} for an outstanding introduction to both topics.  Readers looking for a gentle introduction to mis\`ere quotients may consult the unpublished lecture notes~\cite{siegel_mqlectures}, which include most of the necessary background.

Section~\ref{section:miserequotients} of this paper is an informal exposition of the mis\`ere quotient construction.  In Sections~\ref{section:bipartitemonoids} and~\ref{section:structuretheory}, we formalize the construction and develop a rudimentary structure theory.  Section~\ref{section:mexfunction} dives a bit more deeply into the general structure of mis\`ere quotients, and in Section~\ref{section:nkh} we discuss a connection between the combinatorial structure of normal and mis\`ere play.

Finally, Appendix~\ref{appendix:phivalues} summarizes the solutions to various octal games, obtained using the techniques of this paper.  Several interesting counterexamples are also presented there.  A supplementary appendix, available online on the arXiv~\cite{plambeck_2007supp}, contains further details regarding these solutions as well as descriptions of the algorithms used to calculate them.

\section{Mis\`ere Quotients}

\label{section:miserequotients}

Let $G$ be an impartial game.  We denote by $o^+(G)$ the normal-play outcome of~$G$.  Thus $o^+(G) = \SP$ if second player (the \emph{p}revious player) can force a win, assuming normal play; otherwise, $o^+(G) = \SN$, indicating that first player (the \emph{n}ext player) can force a win.  Likewise, we denote by $o^-(G)$ the mis\`ere-play outcome of $G$.  We say that $G$ is a \emph{normal \PPos} if $o^+(G) = \SP$, etc.

Now let $G$ and $H$ be impartial games.  We define
\[G =^+ H \textrm{ if } o^+(G + X) = o^+(H + X), \textrm{ for every impartial game } X.\]
\[G =^- H \textrm{ if } o^-(G + X) = o^-(H + X), \textrm{ for every impartial game } X.\]
It is customary to write simply $G = H$ instead of $G =^- H$ when we are firmly in the context of mis\`ere games.

A major goal of combinatorial game theory is to understand the structure of the equivalence classes of games modulo equality.  As mentioned in the introduction, Sprague and Grundy completed this program convincingly for normal play.

\begin{fact}[Sprague--Grundy Theorem]
Let $G$ be any impartial game.  Then for some Nim-heap $\Star n$, we have $G =^+ \Star n$.
\end{fact}

The Sprague--Grundy Theorem yields a remarkably simple structure for the normal-play equivalence classes of impartial games.  In mis\`ere play, however, the situation is vastly more complex.  For example, consider just those impartial games born by day 6.  In normal play, there are precisely seven of them, up to equivalence: $0,\Star,\Star 2,\ldots,\Star 6$.  By contrast, Conway has shown that in mis\`ere play there are more than $2^{4171779}$~\cite{conway_1976}.

The solution to this problem is to localize the definition of equality.  Let $\A$ be some fixed set of games, closed under addition, and suppose that whenever $G \in \A$, then every option of $G$ is also in $\A$.  (Typically, $\A$ will be the set of all positions that arise in some fixed game $\Gamma$, such as Dawson's Chess.)  Then define, for all $G,H \in \A$,
\[G \equiv_\A H \textrm{ if } o^-(G + X) = o^-(H + X), \textrm{ for every } X \in \A.\]
It is easily seen that $\equiv_\A$ is an equivalence relation.  It is furthermore a congruence: if $G \equiv_\A H$ and $K \in \A$, then since $\A$ is closed under sums, we have $G + K \equiv_\A H + K$.  Consequently, addition modulo $\equiv_\A$ defines a monoid $\Q$:
\[\Q = \A / \equiv_\A.\]
Moreover, if $G \equiv_\A H$, then necessarily $o^-(G) = o^-(H)$ (taking $X = 0$ in the definition of equivalence).  Therefore each $x \in \Q$ has a well-defined outcome, and we can put
\[\P = \{[G]_{\equiv_\A} : G \in \A,\ o^-(G) = \SP\}.\]
The structure $\QP$ is the \emph{mis\`ere quotient} of $\A$, and we denote it by $\Q(\A)$.  We will refer to $\P$ as the \emph{$\SP$-portion} of $\Q(\A)$.

As we will see, the mis\`ere quotient of $\A$ is often finite, even when $\A$ is infinite.  In such cases, our reduction makes the theory tractable, as follows.  Let $\Phi : \A \to \Q$ be the quotient map, defined by $\Phi(G) = [G]_{\equiv_\A}$.  Suppose we can identify a set of generators $\mathscr{H} \subset \A$, together with the images $\Phi(H) \in \Q$ of each generator $H \in \mathscr{H}$.  (For example, $\A$ might be the set of all Dawson's Chess positions, and $\mathscr{H}$ the set of Dawson's Chess positions consisting of a single contiguous row of pawns.)  Given a position $G \in \A$, we can write~$G$ as a sum of generators, $G = H_1 + H_2 + \cdots + H_k$, and compute the images $x_1 = \Phi(H_1)$, $x_2 = \Phi(H_2)$, $\ldots$, $x_k = \Phi(H_k)$.  Then we simply check whether
\[x_1\cdot x_2\cdot \cdots \cdot x_k \in \P,\]
and this determines $o^-(G)$.  The problem of finding the outcome of the (possibly quite complicated) sum $G$ is therefore reduced to a small number of operations on the finite multiplication table $\Q$.

In this paper, we are principally concerned with \emph{octal games}~\cite{berlekamp_1982,guy_1956}, in which the positions are sums of heaps of tokens.  If $\Gamma$ is an octal game, we will usually denote by $H_n$ a $\Gamma$-heap of size $n$.  We denote by $\A$ the set of all positions of $\Gamma$ (the ``heap algebra'') and by $\mathscr{H}$ the set of single-heap positions.  When convenient, we will regard $\A$ as a free commutative monoid on the countable generating set $\mathscr{H}$.  Finally, it will be useful to write $\mathscr{H}_n = \{H_1,\ldots,H_n\}$ and $\A_n = \<\mathscr{H}_n\>$, the submonoid of $\A$ generated by $\mathscr{H}_n$.  We refer to $\Q(\A_n)$ as the \emph{$n^\textrm{th}$ partial quotient} of $\Gamma$.  For convenience, we will sometimes write $\Q(\Gamma)$ and $\Q_n(\Gamma)$ in place of $\Q(\A)$ and $\Q(\A_n)$, respectively.

To present the solution to $\Gamma$, it suffices to specify three things: the monoid~$\Q$, the $\SP$-portion $\P$, and the quotient map $\Phi : \A \to \Q$.  Since $\Phi$ is a monoid homomorphism, we need only specify the values $\Phi(H_n)$ for each $n$.  We can use this information to reconstruct perfect play for $\Gamma$, as described above.

In this paper we will describe several such solutions to various octal games.  These solutions were obtained with the help of \emph{MisereSolver} \cite{miseresolver}, our software for computing mis\`ere quotients, which is described in detail in the supplement~\cite{plambeck_2007supp}.  The one missing ingredient is the Guy--Smith--Plambeck \emph{Periodicity Theorem}, which enables one to compute $(\Q,\P,\Phi)$ in terms of some finite partial quotient $\Q(\A_n)$.  For the record, we state the theorem here.

\begin{fact}[Periodicity Theorem \cite{guy_1956,plambeck_2005}]
\label{fact:periodicitytheorem}
Let $\Gamma$ be an octal game, let $d$ be the index of its last non-zero code digit, and let $\A$ be the free commutative monoid on the set of heaps $\mathscr{H} = \{H_1,H_2,H_3,\ldots\}$.  Let $\QP$ be the mis\`ere quotient of $\A$, with quotient map $\Phi : \A \to \Q$.  Suppose that
\[\Phi(H_{n+p}) = \Phi(H_n), \textrm{ for all $n$ such that } n_0 \leq n < 2n_0 + p + d.\]
Then
\[\Phi(H_{n+p}) = \Phi(H_n) \textrm{ for all } n \geq n_0.\]
\end{fact}

Therefore, if we compute a partial quotient $\Q(\A_{n_0+2p+d})$ that satisfies the hypotheses of the Periodicity Theorem, then we obtain a complete solution for~$\Gamma$, with $\Q(\A) \cong \Q(\A_{n_0+2p+d})$.

\subsection*{Example: Mis\`ere Guiles}
\suppressfloats[t]

\emph{Guiles} is the octal game \textbf{0.15}.  It was previously unsolved (despite non-trivial effort), but \emph{MisereSolver} can nonetheless dispense a complete solution in under one second.  A presentation for the mis\`ere quotient $\QP$ is shown in Figure~\ref{figure:015}, together with the single-heap values of the quotient map $\Phi$.  The monoid $\Q$ has~42 elements; 12 of these are in $\P$.

\begin{figure}
\centering
\xpretend{10}{a a 1 a a b b a b b a a 1 c c b b d b e c c f c c b g d h i ab2 abg f abg abe b3 h d h h ab2 abe f2 abg abg b3 h d h h ab2 abg f2 abg abg b3 b3 d b3 b3 ab2 abg f2 abg abg b3 b3 d b3 b3 ab2 ab2 f2 ab2 ab2 b3 b3 d b3 b3 ab2 ab2 f2 ab2 ab2}

\[\begin{array}{@{}r@{~}l@{}}
\Q \cong & \present[t]{9.5cm}{a,b,c,d,e,f,g,h,i}{a2=1,b4=b2,bc=ab3,c2=b2,b2d=d,cd=ad,d3=ad2,b2e=b3,de=bd,be2=ace,ce2=abe,e4=e2,bf=b3,df=d,ef=ace,cf2=cf,f3=f2,b2g=b3,cg=ab3,dg=bd,eg=be,fg=b3,g2=bg,bh=bg,ch=ab3,dh=bd,eh=bg,fh=b3,gh=bg,h2=b2,bi=bg,ci=ab3,di=bd,ei=be,fi=b3,gi=bg,hi=b2,i2=b2} \vspace{0.25cm} \\
\P = & \plist[t]{9.5cm}{a,b2,bd,d2,ae,ae2,ae3,af,af2,ag,ah,ai}
\end{array}\]

\caption{The mis\`ere quotient of Guiles (\textbf{0.15}).}
\label{figure:015}
\end{figure}


\subsection*{The Tame Quotients}
\label{tameorwild}

In a sense, the simplest mis\`ere quotients are those that arise in Nim.  Put $\mathcal{T}_0 = \Q(0)$, and for each $n \geq 1$ define
\[\mathcal{T}_n = \Q(\Star 2^{n-1}).\]
Finally, put
\[\mathcal{T}_\infty = \Q(0,\Star,\Star2,\Star3,\ldots),\]
the quotient of mis\`ere Nim.

The quotients $\mathcal{T}_n$ are extremely common, and they correlate with an important property of mis\`ere games known as \emph{tameness}.  Roughly speaking, an impartial game is said to be \emph{tame} if its mis\`ere-play strategy can be described in terms of the mis\`ere-play strategy of Nim.  There are several ways to formalize the definition; Conway~\cite{berlekamp_1982,conway_1976} uses a modification of Grundy values known as \emph{genus symbols}.  For the purposes of this paper, we simply define a set of games~$\A$ to be \emph{tame} if its mis\`ere quotient is isomorphic to $\mathcal{T}_n$ for some~$n$.\footnote{The two definitions are equivalent: $\A$ is tame (in the quotient sense) if and only if every $G \in \A$ is tame (in Conway's sense).  See~\cite{siegel_200Xe} for a proof of this fact, and~\cite{plambeck_2005,plambeck_200X} for further discussion of the relationship between mis\`ere quotients and genus symbols.

Note also that the definition of ``tame'' in \emph{Winning Ways}~\cite{berlekamp_1982} is somewhat more general than Conway's original definition~\cite{conway_1976} (``hereditarily tame'').  However, our hereditary closure assumption on~$\A$ guarantees that the two definitions coincide for our purposes.}  Otherwise, we say that $\A$ is \emph{wild}.


It is worthwhile to write down explicit presentations for the $\mathcal{T}_n$:
\[\begin{array}{l@{~\cong~}l}
\mathcal{T}_1 & \<a~|~a^2=1\> \bigstrut \\
\mathcal{T}_2 & \<a,b~|~a^2=1,\ b^3=b\> \bigstrut \\
\mathcal{T}_n &
\begin{array}[t]{@{}c@{~|~}p{7cm}@{}}
\<a,b_1,b_2,\ldots,b_{n-1} & \raggedright
${a^2=1},\ {b_1^3=b_1},\ {b_2^3=b_2},\ \ldots,\ {b_{n-1}^3=b_{n-1}},\ \allowbreak {b_1^2=b_2^2=b_3^2=\cdots=b_{n-1}^2}\>$
\end{array} \bigstrut \vspace{0.05cm} \\
\mathcal{T}_\infty &
\begin{array}[t]{@{}c@{~|~}p{7cm}@{}}
\<a,b_1,b_2,\ldots & \raggedright
$a^2=1,\ {b_1^3=b_1},\ {b_2^3=b_2},\ \ldots,\ \allowbreak
{b_1^2=b_2^2=b_3^2=\cdots}\>$
\end{array}
\end{array}\]

Now $\mathcal{T}_1 \cong \mathbb{Z}_2$ (as a monoid).  For $n \geq 2$, the structure of $\mathcal{T}_n$ is best described as follows.  Note that $b_1^2 = b_2^2 = \cdots = b_{n-1}^2$; let $z$ be this value.  Then every $x \in \mathcal{T}_n$ satisfies $zx = x$, \emph{except} for $x = 1,a$.  If we let
\[\mathcal{K}_n = \{x \in \mathcal{T}_n : zx = x\},\]
then $\mathcal{K}_n$ is a group with identity $z$, and we have $\mathcal{T}_n = \mathcal{K}_n \cup \{1,a\}$.  Furthermore, $\mathcal{K}_n \cong \mathbb{Z}_2^n$, with generators $za,b_1,b_2,\ldots,b_{n-1}$.  The remaining two elements $\{1,a\}$ form a separate copy of $\mathbb{Z}_2$.  Therefore $|\mathcal{T}_n| = 2^n+2$.

Likewise, we can define $\mathcal{K}_\infty$ in the same way, and we have
\[\mathcal{K}_\infty \cong \bigoplus_{\mathbb{N}} \mathbb{Z}_2,\]
and $\mathcal{T}_\infty = \mathcal{K}_\infty \cup \{1,a\}$.

Readers familiar with Conway's genus theory will recognize $1$ and $a$ as the \emph{fickle units}, and $z$ and $az$ as the \emph{firm units}.

$\mathcal{T}_2$, in particular, is especially common; for example, it is the quotient of {\bf 0.23}, {\bf 0.31} (``Stalking'') and {\bf 0.52}.

\subsection*{The Smallest Wild Quotient}
\label{section:wildquotient}

The smallest wild mis\`ere quotient is an eight-element monoid that we denote by ${\mathcal R}_8$:
\[\mathcal{R}_{8} = \<a,b,c~|~a^2=1, \ b^3=b,\ bc=ab,\ c^2=b^2 \>; \quad \mathcal{P} = \{a,b^2\}.\]
In Section \ref{section:structuretheory} of this paper, we will prove that $\mathcal{T}_0$, $\mathcal{T}_1$ and $\mathcal{T}_2$ are the only quotients smaller than $\mathcal{R}_8$.  It is also possible to show that $\mathcal{R}_8$ is the \emph{unique} mis\`ere quotient of order 8 (up to isomorphism), though this is more difficult; see~\cite{siegel_200Xe}.

Many octal games have mis\`ere quotient $\mathcal{R}_8$, including {\bf 0.75} (Section~\ref{allemang}), {\bf 0.512}, and {\bf 4.56}.



\section{Bipartite Monoids}

\label{section:bipartitemonoids}

In this section we introduce an algebraic framework for the mis\`ere quotient construction.

\begin{definition}
A \emph{bipartite monoid} is a pair $\QP$, where $\Q$ is a monoid and $\P \subset \Q$ is an identified subset of $\Q$.
\end{definition}

For example, let $\A$ be a set of impartial games, closed under addition.  Let $\B \subset \A$ be the set of normal-play $\SP$-positions in $\A$.  Then $(\A,\B)$ is a bipartite monoid.  If instead we take $\B$ to be the set of mis\`ere-play $\SP$-positions, we get a different bipartite monoid.

We will be largely concerned with quotients of bipartite monoids.  The discussion can be framed categorically, as follows.

\begin{definition}
Let $\QP$, $(\Q',\P')$ be bipartite monoids.  A \emph{homomorphism} $f : \QP \to (\Q',\P')$ is a monoid homomorphism $\Q \to \Q'$ such that
\[x \in \P \Longleftrightarrow f(x) \in \P'\]
for all $x \in \Q$.
\end{definition}

Then $(\Q',\P')$ is a quotient of $\QP$ if there exists a surjective homomorphism $\QP \to (\Q',\P')$.

\begin{definition}
Let $\QP$ be a bipartite monoid.  Two elements $x,y \in \Q$ are \emph{indistinguishable} if, for all $z \in \Q$, we have
\[xz \in \P \Longleftrightarrow yz \in \P.\]
\end{definition}

Our motivation for this definition comes from the game-theoretic interpretation.  For example, let $\A$ be the set of all impartial games, and let~$\B$ consist of all normal $\SP$-positions.  Then $X,Y \in \A$ are indistinguishable if and only if $X + Z$ and $Y + Z$ have the same outcomes, for all $Z$.  Thus indistinguishability coincides with Grundy equivalence.

\begin{definition}
\tolerance=300
A bipartite monoid $\QP$ is \emph{reduced} if its elements are pairwise distinguishable.  We write \emph{r.b.m.}\ as shorthand for \emph{reduced bipartite monoid}.
\end{definition}

Every bipartite monoid $\QP$ has a quotient that is reduced, and there is a simple procedure for constructing this r.b.m.  Fix $\QP$, and define, for all $x,y \in \Q$,
\[x~\rho~y \Longleftrightarrow x \textrm{ and } y \textrm{ are indistinguishable}.\]
Now $\rho$ is a congruence: it is an equivalence relation, and if $v~\rho~w$ and $x~\rho~y$, then $vx~\rho~wy$.  Thus the equivalence classes mod $\rho$ form a monoid $\Q'$.  Moreover, if $x~\rho~y$, then $x \in \P \Leftrightarrow y \in \P$ (taking $z = 1$ in the definition of indistinguishability).  Therefore, the set $\P' = \{[x]_\rho : x \in \P\}$ is well-defined, and it is easily seen that $(\Q',\P')$ is a quotient of $\QP$.  Finally, $(\Q',\P')$ must be reduced: if $[x]_\rho$ and $[y]_\rho$ are indistinguishable in $(\Q',\P')$, then $x$ and $y$ are indistinguishable in $(\Q,\P)$, whereupon $[x]_\rho = [y]_\rho$.

\begin{definition}
Let $\QP$ be a bipartite monoid.  The \emph{reduction} of $\QP$ is the r.b.m.\ $(\Q',\P')$ defined by $\Q' = \Q/\rho$, $\P' = \{[x]_\rho : x \in \P\}$.
\end{definition}

\begin{example}
Let $\A$ be a set of impartial games, closed under addition, and let~$\B$ be the set of normal-play $\SP$-positions in $\A$.  Let $\QP$ be the reduction of $(\A,\B)$.
\begin{itemize}
\item If the Grundy values in $\A$ are bounded, say $\{0,1,\ldots,2^n-1\}$, then by the Sprague--Grundy Theorem $\Q \cong \mathbb{Z}_2^n$.
\item If the Grundy values in $\A$ are unbounded, then
\[\Q \cong \bigoplus_{\mathbb{N}} \mathbb{Z}_2.\]
\end{itemize}
In either case, $\P = \{0\}$, and if $\Phi : \A \to \Q$ is the quotient map, then $\Phi(X)$ is the Grundy value of $X$ (in binary).
\end{example}

\begin{example}
Let $\A$ be the set of all impartial games, and let $\B$ be the set of all mis\`ere $\SP$-positions.  Then $X~\rho~Y$ if and only if $X =^- Y$ (if and only if $X$ and~$Y$ have the same mis\`ere canonical form).
\end{example}


We now show that the reduction of $\QP$ is the \emph{unique} reduced quotient of $\QP$, up to isomorphism.

\begin{proposition}
\label{proposition:uniquereduction}
Let $\QP$ and $(\mathcal{S},\mathcal{R})$ be bipartite monoids, with reductions $(\Q',\P')$ and $(\mathcal{S}',\mathcal{R}')$, respectively.  Suppose that $(\mathcal{S},\mathcal{R})$ is a quotient of $\QP$, i.e., there exists a surjective homomorphism $f : \QP \to (\mathcal{S},\mathcal{R})$.  Then there is an isomorphism $i : (\Q',\P') \cong (\mathcal{S}',\mathcal{R}')$ making the following diagram commute:

\[\xymatrix{
\Q \ar@{>>}[r] \ar_{f}[d] & \Q' \ar@{=>}^{i}[d] \\
\mathcal{S} \ar@{>>}[r] & \mathcal{S}'
}\]
\end{proposition}

\begin{proof}
Since $f$ is a surjective homomorphism of bipartite monoids, we have
\begin{eqnarray*}
[x] = [y] & \textrm{iff} & xz \in \P \Leftrightarrow yz \in \P \textrm{ for all } z \in \Q \\
& \textrm{iff} & f(xz) \in \mathcal{R} \Leftrightarrow f(yz) \in \mathcal{R} \textrm{ for all } z \in \Q \\
& \textrm{iff} & f(x)w \in \mathcal{R} \Leftrightarrow f(y)w \in \mathcal{R} \textrm{ for all } w \in \mathcal{S} \\
& \textrm{iff} & [f(x)] = [f(y)].
\end{eqnarray*}
Therefore we can define $i : \Q' \to \mathcal{S}'$ by $i([x]) = [f(x)]$, and the conclusions are apparent.
\end{proof}

\begin{corollary}
A r.b.m.\ has no proper quotients (in the category of bipartite monoids).
\end{corollary}

\section{The Structure of Mis\`ere Quotients}

\label{section:structuretheory}

Let $\A$ be a set of impartial games.  We say that $\A$ is \emph{hereditarily closed} if, whenever $G \in \A$ and $G'$ is an option of $G$, then also $G' \in \A$.  When $\A$ is both closed under addition and hereditarily closed, we simply say that $\A$ is \emph{closed}.  We denote by $\cl(\A)$ the closure of an arbitrary set $\A$ (that is, the smallest closed superset of~$\A$).

\begin{definition}
Let $\A$ be a closed set of impartial games, and let $\B$ be the set of all mis\`ere $\SP$-positions in $\A$.  Then the \emph{mis\`ere quotient of~$\A$}, denoted $\Q(\A)$, is the reduction $\QP$ of $(\A,\B)$.
\end{definition}

For convenience, we will sometimes write $\Q(\A)$ in place of $\Q(\cl(\A))$, even when $\A$ is not closed.  Likewise, if $G$ is a single game, we may write $\Q(G)$ in place of $\Q(\cl(\{G\}))$.

Some aspects of the theory can be generalized to sets that are closed under addition, but are not hereditarily closed.  However, virtually all sets of games that are interesting to us will be hereditarily closed.  For example, if $\A$ is the set of positions that arise in some specified heap game $\Gamma$, then $\A$ is necessarily closed, since all options are to sums of smaller heaps.  Thus there is not much harm in taking closure to be a basic assumption.  Furthermore, giving up closure would require us to loosen the theory considerably.

Notice, for example, that closure implies that $0 \in \A$.  Furthermore, if $\A$ is non-trivial, then it must contain a game of birthday exactly 1.  The only such game is $\Star$, so necessarily $\Star \in \A$.  This gives our first proposition:

\begin{proposition}
Every non-trivial mis\`ere quotient contains an element $a \neq 1$ satisfying $a^2 = 1$.
\end{proposition}

\begin{proof}
Let $a = \Phi(\Star)$.  Since $\Star + \Star = 0$, we have $\Star + \Star \equiv_\A 0$; therefore $a^2 = 1$.  Furthermore, $0$ is an \NPos{} and $\Star$ is a \PPos{}, so $a \neq 1$.
\end{proof}

Now we rattle off some elementary facts about $\Phi$-values:

\begin{proposition}
\label{proposition:ppos}
Let $\QP = \Q(\A)$ be a non-trivial mis\`ere quotient and fix $x \in \Q$.  Then there is some $y \in \Q$ with $xy \in \P$.
\end{proposition}

\begin{proof}
If $x = 1$, then $x\Phi(\Star) = \Phi(\Star) \in \P$.  Otherwise, fix some $G \in \A$ with $x = \Phi(G)$, and consider $G + G$.  If $\Phi(G + G) \in \P$, then we can simply put $y = x$.  If $\Phi(G + G) \not\in \P$, then $G + G$ is an \NPos.  Since $x \neq 1$, we have $G \neq 0$, so there is some $G'$ with $G + G'$ a \PPos.  But then $\Phi(G + G') \in \P$, so we can put $y = \Phi(G')$.
\end{proof}

\begin{corollary}
\label{corollary:optdist}
Let $\QP = \Q(\A)$, let $\Phi : \A \to \Q$ be the quotient map, and fix $G \in \A$.  If $G'$ is an option of $G$, then $\Phi(G') \neq \Phi(G)$.
\end{corollary}

\begin{proof}
Since $G$ has an option, $\A$ must be nontrivial, so by Proposition~\ref{proposition:ppos} there is some $y \in \Q$ with $\Phi(G)y \in \P$.  Fix $H$ with $\Phi(H) = y$, so that $G + H$ is a \PPos.  Then necessarily $G' + H$ is an \NPos, so $\Phi(G')y \not\in \P$, and $y$ distinguishes $\Phi(G')$ from $\Phi(G)$.
\end{proof}

We now establish some interesting results regarding the order of a mis\`ere quotient~$\Q$.  We first show that, except for the trivial quotient $\Q(0)$, every mis\`ere quotient has even order.

\begin{theorem}
Let $\QP$ be a mis\`ere quotient.  If $\QP$ is finite and non-trivial, then~$|\Q|$ is even.
\end{theorem}

\begin{proof}
Write $\QP = \Q(\A)$ for some closed set of games $\A$, and let ${\Phi : \A \to \Q}$ be the quotient map.  Since $\A$ is non-trivial, we have $\Star \in \A$.  Put $a = \Phi(\Star)$, and define $f : \Q \to \Q$ by $f(x) = ax$.

Now for every $G \in \A$, we know that $G + \Star \in \A$.  Since $G$ is an option of $G + \Star$, Corollary~\ref{corollary:optdist} gives $\Phi(G) \neq \Phi(G + \Star)$.  But $\Phi(G + \Star) = a\Phi(G)$, so we conclude that $f(x) \neq x$ for all $x \in \Q$.  Furthermore, since $\Star + \Star = 0$, we have $f(f(x)) = x$ for all $x \in \Q$.  Therefore $f$ induces a perfect pairing of elements of~$\Q$.
\end{proof}

Next we show that there is no mis\`ere quotient of order 4, and that the quotients of orders 1, 2 and 6 are unique (up to isomorphism).

\begin{lemma}
\label{lemma:subalgebras}
Let $\A$ be a closed set of games, let $\QP = \Q(\A)$, and let $\B \subset \A$ be a closed subset.  Then there is a submonoid $\mathcal{R}$ of $\Q$ such that
\[\Q(\B) \textrm{ is the reduction of } (\mathcal{R},\P \cap \mathcal{R}).\]
\end{lemma}

\begin{proof}
Let $\Phi : \A \to \Q$ be the quotient map and put $\mathcal{R} = \{\Phi(G) : G \in \mathscr{B}\}$.  Since $\mathscr{B}$ is closed, $\mathcal{R}$ is a submonoid of $\Q$.  Furthermore, the restriction
\[\Phi \upharpoonright \mathscr{B} : \mathscr{B} \to \mathcal{R}\]
is a surjective homomorphism of bipartite monoids.  Therefore the reduction of $(\mathcal{R},\P \cap \mathcal{R})$ is also a quotient of $\mathscr{B}$.  By uniqueness (Proposition \ref{proposition:uniquereduction}), it must be the mis\`ere quotient of $\mathscr{B}$.
\end{proof}

\begin{lemma}
\label{lemma:canonicaliso}
Let $\A,\B$ be closed sets of games.  Suppose that every $G \in \A$ is canonically equal to some $H \in \B$, and vice versa.  Then $\Q(\A) \cong \Q(\B)$.
\end{lemma}

\begin{proof}
Let $\mathcal{S}$ be the monoid of canonical equivalence classes of mis\`ere games, and let $\mathcal{R} \subset \mathcal{S}$ be the set of canonical \PPos{}s.  (Thus $(\mathcal{S},\mathcal{R})$ is the reduction of the universe of mis\`ere games.)  By assumption on~$\A$ and~$\B$, the natural homomorphisms $\A \to \mathcal{S}$ and $\B \to \mathcal{S}$ have the same image.  Let $\Q \subset \mathcal{S}$ be this image and let $\P = \Q \cap \mathcal{R}$.  Let $(\Q',\P')$ be the reduction of $\QP$.  Then $(\Q',\P')$ is a quotient of both $\A$ and $\B$, so by Proposition~\ref{proposition:uniquereduction}, we have
\[\Q(\A) \cong (\Q',\P') \cong \Q(\B). \qedhere\]
\end{proof}

\begin{lemma}
\label{lemma:star2}
Let $\A$ be a non-empty closed set of games.  Then either:
\begin{enumerate}
\item[(i)] Every $G \in \A$ satisfies either $G = 0$ or $G = \Star$; or
\item[(ii)] There exists a $G \in \A$ with $G = \Star 2$.
\end{enumerate}
\end{lemma}

\begin{proof}
Suppose (i) fails.  Then there is some $G \in \A$ with $G \neq 0,\Star$.  Choose such $G$ with minimal birthday.  Then every $G'$ satisfies either $G' = 0$ or $G' = \Star$.  Thus either $G = \{0\}$, or $G = \{\Star\}$, or $G = \{0,\Star\}$.  But the first two possibilities imply $G = \Star$ and $G = 0$, respectively; by assumption on $G$, neither is true, so $G = \{0,\Star\} = \Star 2$.
\end{proof}

\begin{theorem}
\label{theorem:t2}
Let $\A$ be a non-empty closed set of games.  Then either:
\begin{enumerate}
\item[(i)] $\Q(\A) \cong \mathcal{T}_0$; or
\item[(ii)] $\Q(\A) \cong \mathcal{T}_1$; or
\item[(iii)] There exists a closed subset $\B \subset \A$ such that $\Q(\B) \cong \mathcal{T}_2$.
\end{enumerate}
\end{theorem}

\begin{proof}
\emph{Case 1}: Every $G \in \A$ satisfies either $G = 0$ or $G = \Star$.  If $\A = \{0\}$, then $\Q(\A) = \mathcal{T}_0$, by definition.  Otherwise, $\A$ must contain a game of birthday exactly 1; $\Star$ is the only such game, so $\Star \in \A$.  By Lemma~\ref{lemma:canonicaliso}, $\Q(\A) \cong \Q(\Star) =~\mathcal{T}_1$.

\vspace{0.15cm}\noindent
\emph{Case 2}: Otherwise, by Lemma \ref{lemma:star2}, there is some $G \in \A$ with $G = \Star 2$.  Choose such~$G$ with minimal birthday and let $\B = \cl(\{G\})$.  By minimality of~$G$, every element of $\B$ is equal to some element of $\cl(\{\Star,\Star2\})$.  By Lemma~\ref{lemma:canonicaliso}, $\Q(\B) \cong \Q(\Star,\Star2) = \mathcal{T}_2$.
%
\end{proof}

\begin{corollary}
There is no mis\`ere quotient of order 4, and exactly one each of orders 1, 2 and 6 (up to isomorphism).
\end{corollary}

\begin{proof}
Let $\QP = \Q(\A)$ be a mis\`ere quotient.  Consider each possibility in Theorem \ref{theorem:t2}.  In case (i), $|\Q| = 1$.  In case (ii), $|\Q| = 2$.  Finally, in case (iii), there is some $\mathscr{B} \subset \A$ with $\Q(\mathscr{B}) \cong \mathcal{T}_2$.  By Lemma \ref{lemma:subalgebras}, there is some $\mathcal{R} < \Q$ whose reduction is $\mathcal{T}_2$.  Therefore
\[6 = |\mathcal{T}_2| \leq |\mathcal{R}| \leq |\Q|.\]
Furthermore, if $|\Q| = 6$, then both inequalities collapse and
\[\mathcal{T}_2 \cong (\mathcal{R},\P \cap \mathcal{R}) = (\Q,\P). \qedhere\]
\end{proof}

In fact, there is also just one quotient of order 8, but this is harder to establish.  See \cite{siegel_200Xe} for a proof.









\section{Transition Algebras and the Mex Function}

\label{section:mexfunction}

Let $\A$ be a set of games with quotient $\QP = \Q(\A)$, and let $G$ be a game with $\opts(G) \subset \A$.  Consider the set $\Phi''G \subset \Q$ defined by
\[\Phi''G = \{\Phi(G') : G' \textrm{ is an option of } G\}.\]
What can we say about the extension $\Q(\A \cup \{G\})$?  In normal play, we know that $\Q(\A \cup \{G\}) \cong \Q(\A)$ if and only if $\Phi''G$ excludes at least one Grundy value in $\A$.  Furthermore, in that case we necessarily have $\Phi(G) = \mathrm{mex}(\Phi''G)$.

In this section we investigate analogous questions in mis\`ere play: Given just $\Phi''G$, can we determine whether $\Q(\A \cup \{G\}) \cong \Q(\A)$?  If so, can we determine $\Phi(G)$, again given just $\Phi''G$?

The answer to both of these questions is yes, but we need more information than is contained in the quotient $\Q(\A)$.  The primary goal of this section is to introduce an intermediate structure $T = T(\A)$---the \emph{transition algebra of~$\A$}---that carries exactly the right information to answer these questions.  As we will see, there exists a partial function $F : \Pow(\Q) \to \Q$, depending only on $T$, such that whenever $\opts(G) \subset \A$:
\begin{enumerate}
\item[(i)] $F(\Phi''G)$ is defined iff $\Q(\A \cup \{G\}) \cong \Q(\A)$; and
\item[(ii)] When $F(\Phi''G)$ is defined, then $\Phi(G) = F(\Phi''G)$.
\end{enumerate}
We call $F$ the \emph{mex function for $T$}.

Transition algebras appear to be necessary: there exist sets of games $\A$ and $\B$, with $\Q(\A) \cong \Q(\B)$, whose mex functions are nonisomorphic.  Further, $T$ is finite whenever $\Q$ is finite, so transition algebras indeed provide a useful simplification (as compared to working directly with $\A$).  They have several other applications as well: in attacking the \emph{classification problem}\footnote{That is, the problem of determining the
number of nonisomorphic quotients of order $n$.}~\cite{siegel_200Xe}, and in the design of algorithms for computing mis\`ere quotients~\cite{plambeck_2007supp}.

For the rest of this section, fix a closed set of games $\A$ with mis\`ere quotient $\QP = \Q(\A)$ and quotient map $\Phi : \A \to \Q$.  We begin with an auxiliary definition.

\begin{definition}
Let $x \in \Q$.  Then the \emph{meximal set of $x$ in $\Q$} is given by
\[\mathcal{M}_x = \{y \in \Q : \textrm{there is no } z \in \Q \textrm{ such that both } xz,yz \in \P\}.\]
\end{definition}

\begin{example}
Let $\Q = \mathcal{T}_2 = \<a,b~|~a^2=1,\ b^3=b\>$, let $\P = \{a,b^2\}$, and let $x = b^2$.  Then $\mathcal{M}_x = \{b,ab,ab^2\}$.
Let $G = \Star 2 + \Star 2$ and $H = \Star 2 + \Star 2 + \Star + \Star$.  We have $\Phi(G) = \Phi(H) = x$ (and in fact, $G = H$ canonically).  However, $\Phi''G = \{b,ab\} \subsetneq \mathcal{M}_x$, while $\Phi''H = \{b,ab,ab^2\} = \mathcal{M}_x$.
\end{example}

Meximal sets are motivated by the following lemma.

\begin{lemma}
\label{lemma:meximal}
Let $G \in \A$ and put $x = \Phi(G)$.  Then $\Phi''G \subset \mathcal{M}_x$.
\end{lemma}

\begin{proof}
Suppose (for contradiction) that there is some $y \in \Phi''G$ with $y \not\in \mathcal{M}_x$.  By definition of $\mathcal{M}_x$, there must be some $z$ such that both $xz,yz \in \P$.  Fix $G' \in \opts(G)$ and $X \in \A$ such that $\Phi(G') = y$ and $\Phi(X) = z$.  Then $G + X$ and $G' + X$ are both \PPos s, a contradiction.
\end{proof}


We will soon show a converse: if $\Phi''G = \mathcal{M}_x$, then necessarily $\Phi(G) = x$.
First we introduce the transition algebra of $\A$.

\begin{definition}
The \emph{transition algebra of $\A$} is the set of pairs
\[T(\A) = \{(\Phi(G),\Phi''G) : G \in \A\}.\]
We define the map $\Psi : \A \to T(\A)$ by $\Psi(G) = (\Phi(G),\Phi''G)$.
\end{definition}

\begin{remark}
Since $T \subset \Q \times \Pow(\Q)$, we have that $|T| \leq |\Q| \cdot 2^{|\Q|}$.  In particular, if~$\Q$ is finite, then so is $T$.
\end{remark}

The image $\Psi(G)$ of a game $G$ identifies not just $\Phi(G)$, but also the $\Phi$-values of all options of~$G$.  Thus $\Psi(G)$ determines all possible \emph{transitions} from $G$ to its options, as projected down to the mis\`ere quotient.  The following lemma shows that $T$ has a commutative monoid structure with identity $(1,\emptyset)$.

\begin{lemma}
If $(x,\mathcal{E}),(y,\mathcal{F}) \in T$, then so is $(xy,y\mathcal{E} \cup x\mathcal{F})$.
\end{lemma}

\begin{proof}
If $\Psi(G) = (x,\mathcal{E})$ and $\Psi(H) = (y,\mathcal{F})$, then $\Psi(G + H) = (xy,y\mathcal{E} \cup x\mathcal{F})$.
\end{proof}

The projection map $T \to \Q$ given by $(x,\mathcal{E}) \mapsto x$ is a monoid homomorphism, and this makes $T$ into a bipartite monoid with reduction $\Q$.  Further, it is trivially verified that $\Psi$ is a homomorphism of bipartite monoids, and therefore (by Proposition~\ref{proposition:uniquereduction}) the following diagram commutes:
\[\xymatrix{
\A \ar^<<<<<{\Psi}@{->>}[r] \ar_{\Phi}@{->>}[dr] & T(\A) \ar@{->>}[d] \\ & \Q(\A)
}\]

Transition algebras have several important uses.  The first is given by the following lemma, which establishes the existence of the mex function $F$ and shows that it depends only on $T$.

\begin{theorem}[Generalized Mex Rule]
\label{theorem:generalizedmex}
Let $T = T(\A)$ and let $G$ be a nonzero game with $\opts(G) \subset \A$.  The following are equivalent, for $x \in \Q$.
\begin{enumerate}
\item[(a)] $\Q(\A \cup \{G\}) \cong \Q(\A)$ and $\Phi(G) = x$.
\item[(b)] The following two conditions are satisfied:
\begin{enumerate}
\item[(i)] $\Phi''G \subset \mathcal{M}_x$; and
\item[(ii)] For each $(y,\mathcal{E}) \in T$ and each $n \geq 0$ such that $x^{n+1}y \not\in \P$, we have either: $x^{n+1}y' \in \P$ for some $y' \in \mathcal{E}$; or else $x^nx'y \in \P$ for some $x' \in \Phi''G$.
\end{enumerate}
\end{enumerate}
\end{theorem}

\begin{proof}
(a) $\Rightarrow$ (b): (i) is just Lemma~\ref{lemma:meximal}.  For (ii), fix $(y,\mathcal{E}) \in T$ and $n \geq 0$, and suppose $x^{n+1}y \not\in \P$.  Fix $Y \in \A$ such that $\Psi(Y) = (y,\mathcal{E})$.  Then $(n+1) \cdot G + Y$ is an \NPos, so either $(n+1) \cdot G + Y'$ is a \PPos, or else $n \cdot G + G' + Y$ is a \PPos.  But these imply, in turn, that $x^{n+1}\Phi(Y') \in \P$ and $x^n\Phi(G')y \in \P$.




\vspace{0.15cm}\noindent
(b) $\Rightarrow$ (a): Fix $H \in \A$ with $\Phi(H) = x$, and write $\B = \cl(\A \cup \{G\})$.  We will show that
\[o^-(G + X) = o^-(H + X), \textrm{ for all } X \in \B.  \tag{\dag} \label{eq:verequiv}\]
This implies that $G \equiv_\B H$, so that $\Q(\B) \cong \Q(\A)$ and $\Phi(G) = x$.

Assume (for contradiction) that (\ref{eq:verequiv}) fails, and fix $K \in \B$ with $o^-(G+K) \neq o^-(H+K)$.  Since $\B = \cl(\A \cup \{G\})$, we can write $K = n \cdot G + Y$, for some $n \geq 0$ and $Y \in \A$.  Choose such $K$ with $n$ minimal; and having done so, choose~$K$ to minimize the birthday of~$Y$.  In particular, this implies that no option of~$K$ witnesses the failure of (\ref{eq:verequiv}).  We now have two cases.

\vspace{0.15cm}\noindent
\emph{Case 1}: $G + K$ is an \NPos, but $H + K$ is a \PPos.  Now $G + K'$ cannot be a \PPos, for any option $K'$ of~$K$: for then minimality of $K$ would imply that $H + K'$ is also a \PPos, contradicting the assumption on $H + K$.  Therefore $n = 0$ (so that $K \in \A$), and $G' + K$ is a \PPos{} for some option $G'$ of~$G$.  We conclude that
\[\Phi(H)\Phi(K) = x\Phi(K) \in \P \qquad \textrm{and} \qquad \Phi(G')\Phi(K) \in \P.\]
Thus $\Phi(G') \not\in \mathcal{M}_x$, contradicting condition (b)(i).

\vspace{0.15cm}\noindent
\emph{Case 2}: $G + K$ is a \PPos, but $H + K$ is an \NPos.  Observe that, by (repeated application of) the minimality of $n$, we have
\[o^-(n \cdot G + X) = o^-(n \cdot H + X), \textrm{ for all } X \in \A. \tag{\ddag} \label{eq:n-1equiv}\]
Now $H + n \cdot G + Y = H + K$ is an \NPos{} and $H + Y \in \A$, so by (\ddag) we have that $H + n \cdot H + Y$ is an \NPos.  Put $(y,\mathcal{E}) = \Psi(Y)$.  By assumption, $\Phi(H + n \cdot H + Y) = x^{n+1}y \not\in \P$.  To complete the proof, we show that this contradicts condition (b)(ii).

First, if $y' \in \mathcal{E}$, then $y' = \Phi(Y')$ for some option $Y'$ of~$Y$.  But since $G + K$ is a \PPos, we know that $G + n \cdot G + Y'$ is an \NPos.  By minimality of $K$, this implies $H + n \cdot G + Y'$ is an \NPos, and by (\ddag), $H + n \cdot H + Y'$ is also an \NPos.  Therefore $x^{n+1}y' \not\in \P$.  So $y'$ cannot fulfill condition (b)(ii).

Finally, suppose $x' \in \Phi''G$, and write $x' = \Phi(G')$.  Then $G' + n \cdot G + Y$ is an \NPos.  But $G' + Y \in \A$, so by (\ddag) we know that $G' + n \cdot H + Y$ is also an \NPos.  Therefore $x^nx'y \not\in \P$, so $x'$ cannot fulfill condition (b)(ii).  This completes the proof.
\end{proof}

{\tolerance=300
\begin{corollary}[Mex Interpolation Principle]
\label{corollary:mexinterpolation}
Let $G$ be a game with $\mathrm{opts}(G) \subset \A$, and suppose that
\[\Phi''H \subset \Phi''G \subset \mathcal{M}_x,\]
for some $H \in \A$ ($H$ not identically $0$) with $\Phi(H) = x$.  Then $\Q(\A \cup \{G\}) \cong \Q(\A)$ and $\Phi(G) = x$.
\end{corollary}
}

\begin{proof}
It suffices to show that $G$ satisfies conditions (i) and (ii) in the Generalized Mex Rule.  (i) is assumed.  Now condition (ii) must hold for $H$ (again by the Generalized Mex Rule); since $\Phi''H \subset \Phi''G$, it must hold for $G$ as well.
\end{proof}

\begin{corollary}
Assume $\A$ is nontrivial, and let $G$ be a game with $\opts(G) \subset \A$ and $\Phi''G = \mathcal{M}_x$.  Then necessarily $\Q(\A \cup \{G\}) \cong \Q(\A)$ and $\Phi(G) = x$.
\end{corollary}

\begin{proof}
Fix $H \in \A$ with $\Phi(H) = x$.  Since $\A$ is nontrivial, we may choose an $H$ that is not identically zero: if $x = 1$ then we can choose $H = \Star + \Star$.  By Lemma~\ref{lemma:meximal}, we have $\Phi''H \subset \mathcal{M}_x$.  Now the corollary follows from the Interpolation Principle.
\end{proof}

\begin{figure}
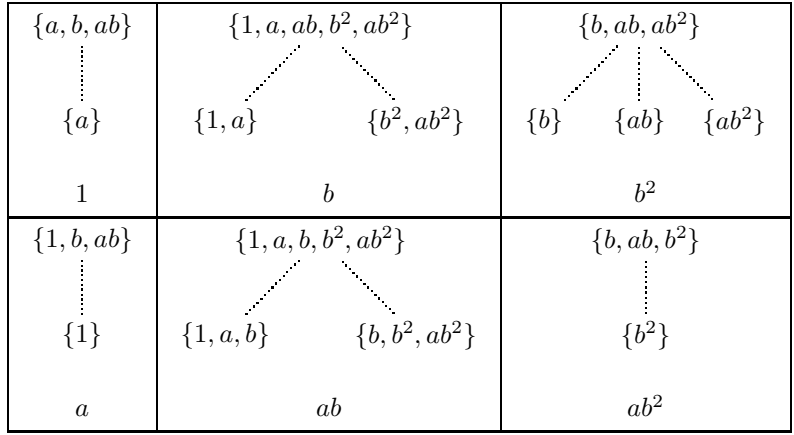

\centering
\begin{tabular}{|c|c|c|}
\hline
  $\xygraph{*+{\{a,b,ab\}} -@{.}[d] *+{\{a\}}}$
& $\xygraph{*+{\{1,a,ab,b^2,ab^2\}} (-@{.}[ld] *+{\{1,a\}}, -@{.}[rd] *+{\{b^2,ab^2\}})}$
& $\xygraph{*+{\{b,ab,ab^2\}} (-@{.}[ld] *+{\{b\}}, -@{.}[d] *+{\{ab\}}, -@{.}[rd] *+{\{ab^2\}})}$
\\ & & \\
$1$ & $b$ & $b^2$ \bigstrut
\\ \hline
  $\xygraph{*+{\{1,b,ab\}} -@{.}[d] *+{\{1\}}}$
& $\xygraph{*+{\{1,a,b,b^2,ab^2\}} (-@{.}[ld] *+{\{1,a,b\}}, -@{.}[rd] *+{\{b,b^2,ab^2\}})}$
& $\xygraph{*+{\{b,ab,b^2\}} -@{.}[d] *+{\{b^2\}}}$
\\ & & \\
$a$ & $ab$ & $ab^2$ \bigstrut
\\ \hline
\end{tabular}
\caption{Schematic of the mex function for $T(\Star 2)$.  Each of the six elements of~$\mathcal{T}_2$ is shown together with its meximal set and antichain of lower bounds.}
\label{figure:mexforT2}
\end{figure}

The Interpolation Principle yields substantial information about the mex function $F$ for $T$.  For a fixed $x \in \Q$, consider the collection
\[\Xi = F^{-1}(x) = \{\mathcal{E} \subset \Q : \mathcal{E} \neq \emptyset \textrm{ and } F(\mathcal{E}) = x\},\]
ordered by inclusion.  (Note that the clause $\mathcal{E} \neq \emptyset$ only matters when $x = 1$, since necessarily $F(\emptyset) = 1$.)  By the Interpolation Principle, we know that $\mathcal{M}_x \in \Xi$, and by Lemma~\ref{lemma:meximal} it serves as an absolute upper bound for $\Xi$.  Now let $\Lambda$ be the antichain of minimal elements of $\Xi$.  Again by the Interpolation Principle, $\Xi$ is completely determined by $\Lambda$ and $\mathcal{M}_x$: indeed,
\[\Xi = \{\mathcal{E} \subset \Q : \mathcal{L} \subset \mathcal{E} \subset \mathcal{M}_x, \textrm{ for some } \mathcal{L} \in \Lambda\}.\]
Therefore $\Xi$ has the structure of an upward-closed subset of the complete Bool\-ean lattice on $\mathcal{M}_x$.  As an example, Figure~\ref{figure:mexforT2} illustrates the structure of $\Xi$ in $T(\Star 2)$ for each of the six elements $x \in \mathcal{T}_2 = \Q(\Star 2)$.

It is a remarkable fact that the upper bound $\mathcal{M}_x$ depends only on $\Q$.  By contrast, the set $\Lambda$ of lower bounds might depend on the fine structure of $T$.

We conclude with a note of caution.  The Generalized Mex Rule asserts only that $\Q(\A \cup \{G\}) \cong \Q(\A)$.  It need \emph{not} be the case that $T(\A \cup \{G\}) = T(\A)$.  Furthermore, the mex function $F^+$ for $T(\A \cup \{G\})$ need not be the same as the mex function $F$ for $T(\A)$: it is possible that $\mathrm{dom}(F^+) \subsetneq \mathrm{dom}(F)$.  However, we \emph{can} conclude that $F = F^+$ in the special case of interpolation:

\begin{proposition}[Mex Interpolation Principle, Strong Form]
\label{proposition:strongmexinterpolation}
Let $G$ be a game with $\mathrm{opts}(G) \subset \A$, and suppose that
\[\Phi''H \subset \Phi''G \subset \mathcal{M}_x,\]
for some $H \in \A$ ($H \neq 0$) with $\Phi(H) = x$.  Let $F$ and $F^+$ be the mex functions for $T(\A)$ and $T(\A \cup \{G\})$, respectively.  Then $F = F^+$.
\end{proposition}

\begin{proof}
We must show that, for each $\mathcal{D} \subset \Q$ and $w \in \Q$,
\[F(\mathcal{D}) = w \quad \textrm{iff} \quad F^+(\mathcal{D}) = w.\]
Now by the Generalized Mex Rule, $F(\mathcal{D}) = w$ iff $\mathcal{D} \subset \mathcal{M}_w$ and
\[\parbox{11cm}{For each $(y,\mathcal{E}) \in T(\A)$ and each $n \geq 0$ such that $w^{n+1}y \not\in \P$, either: $w^{n+1}y' \in \P$ for some $y' \in \mathcal{E}$; or else $w^nw'y \in \P$ for some $w' \in \mathcal{D}$.} \tag{\dag} \]
Likewise, $F^+(\mathcal{D}) = w$ iff $\mathcal{D} \subset \mathcal{M}_w$ and
\[\parbox{11cm}{For each $(y,\mathcal{E}) \in T(\A \cup \{G\})$ and each $n \geq 0$ such that $w^{n+1}y \not\in \P$, either: $w^{n+1}y' \in \P$ for some $y' \in \mathcal{E}$; or else $w^nw'y \in \P$ for some $w' \in \mathcal{D}$.} \tag{\ddag} \]
Now $T(\A) \subset T(\A \cup \{G\})$, so (\ddag) $\Rightarrow$ (\dag).  To prove that (\dag) $\Rightarrow$ (\ddag), it suffices to show the following: whenever $(y,\mathcal{E}) \in T(\A \cup \{G\})$, then there is some $\mathcal{E}' \subset \mathcal{E}$ with $(y,\mathcal{E}') \in T(\A)$.  But this is a simple consequence of the assumptions on $G$: for any $Y \in \A$ and $n \geq 0$, we have $\Phi(n \cdot H + Y) = \Phi(n \cdot G + Y)$, and furthermore $\Phi''(n \cdot H + Y) \subset \Phi''(n \cdot G + Y)$.  This completes the proof.
\end{proof}

Finally, note that if $\Q(\A \cup \{G\}) \neq \Q(\A)$, then the content of $\Q(\A \cup \{G\})$ \emph{cannot} be determined on the basis of $\Phi''G$ alone: it is sensitive to the most minute structural details of $\cl(\A \cup \{G\})$.

\section{The Kernel and Normal Play}

\label{section:nkh}

The strategy for mis\`ere Nim has been known since Bouton~\cite{bouton_1902}, and it is usually formulated as follows:

\begin{center}
\framebox{\parbox{\textwidth-1.5in}{
\centering
Play normal Nim until your move would leave a position consisting entirely of heaps of size 1.  Then play to leave an odd number of heaps of size 1.
}}
\end{center}

In this section we introduce a suitable generalization of this strategy to a wide class of games.  First we recall some basic facts from commutative semigroup theory.  Proofs of these results can be found in a standard reference such as~\cite{grillet_2001}.



Let $\Q$ be a monoid and fix $x,y \in \Q$.  We say that $x$ \emph{divides} $y$, and write $x|y$, if $xz = y$ for some $z \in \Q$.  $x$ and $y$ are \emph{mutually divisible} if $x|y$ and $y|x$.  It is easy to see that mutual divisibility is a congruence.  The congruence class of~$x$ is called the \emph{mutual divisibility class} of $x$.

An \emph{idempotent} is an element $x \in \Q$ such that $x^2 = x$.  If $x$ is an idempotent, then its mutual divisibility class $\mathcal{G}$ is a \emph{group} with identity $x$.  Furthermore, $\mathcal{G}$ is maximal among \emph{group}s contained in~$\Q$.

If $\Q$ is finite, then we can enumerate its idempotents $z_1,z_2,\ldots,z_k$.  The mutual divisibility class $\mathcal{K}$ of their product $z = z_1z_2\cdots z_k$ is called the \emph{kernel} of $\Q$.  Multiplication by $z$ defines a surjective homomorphism $\Q \to \mathcal{K}$, and any surjective homomorphism from $\Q$ onto a group $\mathcal{G}$ factors through this map.  Thus $\mathcal{K}$ is the \emph{group of differences} obtained by adjoining formal inverses to $\Q$.

\begin{proposition}
\label{proposition:kernel_contains_ppos}
Let $\QP = \Q(\A)$ be a finite mis\`ere quotient and let $\mathcal{K}$ be the kernel of $\Q$.  Then $\mathcal{K} \cap \P$ is nonempty.
\end{proposition}

\begin{proof}
Let $z$ be the identity of $\mathcal{K}$.  If $z \in \P$, then we are done.  Otherwise, fix $G$ with $\Phi(G) = z$.  Then $\Phi(G + G) = z^2 = z$.  Now $G + G$ is an \NPos, so some option $G + G'$ must be a \PPos.  But
\[\Phi(G + G') = z\Phi(G') \in \mathcal{K}. \qedhere \]
\end{proof}

\begin{definition}
A finite r.b.m.\ $\QP$ is said to be \emph{regular} if $|\mathcal{K} \cap \P| = 1$, and \emph{normal} if $\mathcal{K} \cap \P = \{z\}$.
\end{definition}

\begin{definition}
The quotient map $\Phi : \A \to \Q$ is said to be \emph{faithful} if
\[\Phi(G) = \Phi(H) \Longrightarrow G \textrm{ and } H \textrm{ have the same normal-play Grundy value}.\]
If $\Phi$ is faithful and $\QP$ is regular (resp.~normal), then we say that $\Phi$ is \emph{faithfully regular} (resp.~\emph{faithfully normal}).
\end{definition}

There do exist irregular quotients; see Appendix~\ref{appendix:phivalues} for an example.  However, they are extremely rare, and in fact most known quotients are normal (including all known full solutions for octal games).  We are not aware of any example of an unfaithful quotient map.  Note that faithfulness depends only on the transition algebra~$T(\A)$, so it is a slightly more robust property than appears at first glance.

\begin{theorem}
\label{theorem:nkh}
Let $\A$ be a closed set of games with finite regular mis\`ere quotient $\QP$ and faithful quotient map $\Phi : \A \to \Q$.  Let $\mathcal{K}$ be the kernel of $\Q$ and let~$z$ be the identity of $\mathcal{K}$.  Then for all $G,H \in \A$, we have
\[z\Phi(G) = z\Phi(H) \Longleftrightarrow G \textrm{ and } H \textrm{ have the same normal-play Grundy value}.\]
\end{theorem}

\begin{proof}
For the forward direction, suppose $z\Phi(G) = z\Phi(H)$.  Fix $Z \in \A$ with $\Phi(Z) = z$. Then $\Phi(G + Z) = \Phi(H + Z)$.  Since $\Phi$ is faithful, we know that $G + Z$ and $H + Z$ have the same Grundy value.  Therefore so do $G$ and $H$.

To complete the proof, we must show that if $G,H \in \A$ have the same Grundy value, then $z\Phi(G) = z\Phi(H)$.

First of all, since $\Q$ is regular, we know that $\mathcal{K} \cap \P$ contains a unique element~$y$.  Fix $Y$ with $\Phi(Y) = y$ and let $r$ be the Grundy value of $Y$.  Suppose $G \in \A$ is any game of Grundy value $r$; we will show that $z\Phi(G) = y$.

Fix $Z \in \A$ with $\Phi(Z) = z$, and put $X = G + Z + Z$.  Since $z$ is an idempotent, we have $\Phi(X) = z\Phi(G)$.  Furthermore, every option of $X$ has the form $G' + Z + Z$ or $G + Z' + Z$, so $\Phi''X \subset \mathcal{K}$.  But since $X$ has Grundy value $r$, none of its options can have Grundy value $r$, so $y \not\in \Phi''X$.  Since $y$ is the unique element of $\mathcal{K} \cap \P$, we infer that $(\Phi''X) \cap \P = \emptyset$, whence $X$ is a \PPos.  Thus $\Phi(X) = z\Phi(G) \in \P$, and by uniqueness of $y$ we conclude that $z\Phi(G) = y$.

Now fix any $G,H \in \A$ and suppose $G$ and $H$ have the same Grundy value.  Then $G + G + Y$ and $G + H + Y$ both have Grundy value $r$, so by the above argument, $z\Phi(G + G + Y) = z\Phi(G + H + Y) = y$.  Put $x = \Phi(G + Y)$; then
\[z\Phi(G + G + Y) = zx\Phi(G) \quad \textrm{and} \quad z\Phi(G + H + Y) = zx\Phi(H).\]
Since $\mathcal{K}$ is a group, there is some $w$ with $wx = z$, and it follows that $z\Phi(G) = z\Phi(H)$, as needed.
\end{proof}

Theorem~\ref{theorem:nkh} yields a one-to-one correspondence between elements of~$\mathcal{K}$ and normal-play Grundy values of games in~$\A$.  It follows that the normal-play structure of~$\A$ is exactly captured by the kernel~$\mathcal{K}$.

For example, in Section~\ref{section:miserequotients} we described the structure of the partial quotients for Nim: if $n \geq 2$, then $\Q(\Star 2^{n-1}) = \{1,a\} \cup \mathcal{K}_n$, where $\mathcal{K}_n \cong \mathbb{Z}_2^n$.  Now if $G \in \cl(\Star 2^{n-1})$, then $\Phi(G) \in \mathcal{K}_n$ if and only if~$G$ contains at least one Nim-heap of size $\geq 2$.  Otherwise, $\Phi(G) = 1$ or~$a$.  Thus positions with at least one heap of size $\geq 2$ map down to normal play, while positions with all heaps of size $\leq 1$ require more delicate consideration; and we recover Bouton's strategy for mis\`ere Nim.  This points the way to the promised generalization, which works for all faithfully normal games~$\Gamma$:

\vspace{0.15cm}

\framebox{\parbox{\textwidth-0.45in}{
\centering
Play normal $\Gamma$ until your move would leave a position outside of $\mathcal{K}$.\\Then pay attention to the fine structure of the mis\`ere quotient.
}}

\vspace{0.15cm}

The difference, of course, is that for mis\`ere Nim this ``fine structure'' is painfully simple, while for arbitrary $\Gamma$ it can be quite complicated indeed.

When $\Gamma$ is faithfully regular but abnormal, Grundy values still carry sufficient information to play optimally inside $\mathcal{K}$, but the actual winning moves are different than in normal play.  For example, if $\mathcal{K} \cap \P = \{az\}$, where $a = \Phi(\Star)$, then the winning move inside $\mathcal{K}$ is to any position of Grundy value~1.

Regularity also has strong consequences for the structure of $\QP$.

\begin{lemma}
\label{lemma:nkhconsequences}
Let $\QP = \Q(\A)$ be a finite regular mis\`ere quotient with faithful quotient map $\Phi : \A \to \Q$.  Let~$\mathcal{K}$ be the kernel of~$\Q$ and let $z$ be the identity of~$\mathcal{K}$.
\begin{enumerate}
\item[(a)] If $x \in \Q$ is an idempotent and $\Phi(G) = x$, then $G$ has Grundy value $0$.
\item[(b)] $\mathcal{K} \cong \mathbb{Z}_2^n$, for some $n$.
\item[(c)] Suppose $\A$ is the set of positions in some heap game $\Gamma$.  If $\Gamma$ is mis\`ere-play periodic, then it is normal-play periodic, and its normal period divides its mis\`ere period.
\end{enumerate}
\end{lemma}

\begin{proof}
(a) $G + G$ necessarily has Grundy value $0$.  But $\Phi(G + G) = x^2 = x$, so
\[z\Phi(G) = zx = z\Phi(G + G),\]
whence by Theorem~\ref{theorem:nkh} $G$ has Grundy value $0$.

\vspace{0.15cm}\noindent
(b) Let $x \in \mathcal{K}$ and fix any $G$ with $\Phi(G) = x$.  Then $G + G$ has Grundy value $0$, so by Theorem~\ref{theorem:nkh} $z\Phi(0) = z\Phi(G+G)$.  Therefore $z = zx^2 = x^2$.  Since $x$ was arbitrary, this shows that every element of $\mathcal{K}$ has order 2.

\vspace{0.15cm}\noindent
(c) Suppose that $\Gamma$ is mis\`ere-play periodic.  Then for some $n_0,p$, we have $\Phi(H_{n+p}) = \Phi(H_n)$ for all $n \geq n_0$.  Therefore $z\Phi(H_{n+p}) = z\Phi(H_n)$ for all $n \geq n_0$.  By Theorem~\ref{theorem:nkh}, $H_{n+p}$ and $H_n$ have the same Grundy value, for all $n \geq n_0$, and the conclusion follows.
\end{proof}

\appendix

\section{The $\Phi$-Values of Various Games}

\label{appendix:phivalues}

We summarize extensive computations obtained using \emph{MisereSolver} \cite{miseresolver}, an extension to the \CGSuite{} computer algebra system \cite{cgsuite}.  The algorithms that drive \emph{MisereSolver} are described in the supplement~\cite{plambeck_2007supp}.

Many of the games presented here were previously unsolved.  Several others had been solved using other methods, but the mis\`ere quotient techniques yield cleaner solutions with much less effort.  We have noted in the text each case in which we are aware of a prior solution.



\subsection{General comments}


\minisection{Many details omitted}
The quotients presented in this paper represent relatively simple examples.  The most complicated quotients we've computed involve thousands of elements, and to write out their minimal presentations would require several largely unenlightening pages.  Such presentations are available on our website~\cite{miseregames_www}, and they're easily reproducible with \emph{MisereSolver}~\cite{miseresolver}, so we've omitted them for most quotients with more than 50 elements.

\minisection{Notation for mis\`ere canonical forms}
Several of the most interesting examples have the form $\Q(G)$, for some specific game~$G$ in mis\`ere canonical form.  To describe such games, we use a slightly modified form of the notation introduced by Conway~\cite[Chapter 12]{conway_1976}.  The notation $\Star n$, where $n$ is a single-digit number, represents (as always) a Nim-heap of size~$n$.  If $G$, $H$, $K$ are mis\`ere canonical forms, stripped of their preceding $\Star$'s, then $\Star GHK$ is the game whose options are to $\Star G$, $\Star H$ and $\Star K$.  Parentheses are used to denote sub-options: $\Star (GH)K$ has options to $\Star GH$ and~$\Star K$.  Finally, $\Star G_\sh$ denotes the game whose only option is to $\Star G$.

For example, $\Star 2_\sh320$ has four options: $0$, $\Star2$, $\Star3$, and the game $\Star 2_\sh$ whose only option is to $\Star 2$.

\minisection{Partial quotients and infinite quotients}
In Section~\ref{section:miserequotients}, we defined the $n^\mathrm{th}$ \emph{partial quotient} $\Q_n(\Gamma)$, obtained by requiring all heaps to have size at most~$n$.  It appears that many octal games have some infinite partial quotients.  Currently, \emph{MisereSolver} is limited to computations on finite mis\`ere quotients; when $\Q_n(\Gamma)$ is infinite, the software goes into an infinite loop at heap size $n$, analyzing successively larger finite approximations to $\Q_n$.

Unfortunately, it is often difficult to prove by hand that a quotient is infinite.  This means that in most cases, we can only speculate whether a quotient is truly infinite, or whether we just haven't run \emph{MisereSolver} for long enough.  (However, see Section~\ref{subsection:infinitequotient} for a proof that $\Q_5(\textbf{0.31011})$ is infinite.)

Figure~\ref{figure:unsolvedtwodigits} gives the largest known partial quotient for every unsolved two-digit octal game.  In each case, we are reasonably sure that the \emph{next} partial quotient is infinite, but a solution to the following problem would give us more confidence.

\begin{problem}
Specify an algorithm to determine whether a given quotient $\Q(\A)$ is infinite.
\end{problem}

A classical theorem of R\'edei~\cite{grillet_2001,redei_1965} implies that every finitely generated commutative monoid is finitely presented.  Since each partial quotient $\Q_n(\Gamma)$ is necessarily finitely generated, there is hope that our techniques can be extended to obtain solutions for many octal games with infinite partial quotients.

\begin{problem}
Specify an algorithm to compute a presentation for the (possibly infinite) quotient $\Q(\A)$, whenever $\A$ is finitely generated.
\end{problem}

\subsection{Wild Octal Games with Known Solutions}
\label{2digits}
We now summarize the wild two- and three-digit octal games whose mis\`ere-play solutions are known.  The results are tabulated in Figure \ref{figure:solvedwildoctals}.  Each row of the table shows a game code (or a schema of equivalent codes), together with the period, preperiod, and quotient order of the mis\`ere solution.  We have also included an appropriate attribution in each case where we are aware of a prior solution.

\begin{figure}
\centering
\begin{tabular}{|>{\bfseries}l|rrrll|}
\hline
\mdseries Code & pd & ppd & $|\Q|$ & Comments & Solved By \\ \hline
0.15 & 10 & 66 & 42 & Guiles, \S \ref{guiles} & \\
0.2\.6\=0 & \multicolumn{2}{c}{\ \dag} & $\infty$ & & Allemang \cite{allemang_2002} \\
0.34\.0 & 8 & 7 & 12 & \S \ref{flanigan}; $\Q \cong \mathcal{S}_{12}$ & Flanigan \cite{berlekamp_1982} \\
0.4\"4\=0 & 24 & 143 & 40 & Duplicate \textbf{0.77} & Sibert \cite{conway_1992} \\
0.53\.0 & 9 & 21 & 16 & \S \ref{allemang} & Allemang \cite{allemang_2002} \\
0.57\.0 & \multicolumn{2}{c}{\ \dag} & $\infty$ & Duplicate \textbf{4.7} & Allemang \cite{allemang_2002} \\
0.71\.0 & 6 & 3 & 36 & \S \ref{allemang} & Flanigan \cite{berlekamp_1982} \\ 
0.72\=0 & 4 & 16 & 24 & \S \ref{allemang} & Allemang \cite{allemang_2002} \\
0.75\=0 & 2 & 8 & 8 & \S \ref{allemang}; $\Q \cong \mathcal{R}_8$ & Allemang \cite{allemang_2002} \\
0.77\.0 & 12 & 71 & 40 & Kayles, \S \ref{kayles} & Sibert \cite{conway_1992} \\
4.\"4\=0 & 12 & 71 & 40 & Cousins of \textbf{0.77} & Sibert \cite{conway_1992} \\
4.7\.0 & \multicolumn{2}{c}{\ \dag} & $\infty$ & Knots, \S \ref{allemang} & Allemang \cite{allemang_2002} \\
0.04\"4 & 24 & 142 & 40 & Duplicate \textbf{0.77} & Sibert \cite{conway_1992} \\
0.07\=4 & 24 & 142 & 40 & Duplicate \textbf{0.77} & Sibert \cite{conway_1992} \\
0.115 & 14 & 92 & 42 & $\frac{7}{5}$-plicate \textbf{0.15} & \\
0.123 & 5 & 5 & 20 & & Plambeck \cite{plambeck_2005} \\
0.144 & 10 & 12 & 30 &  & \\
0.152 & 48 & 25 & 34 &  & \\
0.153 & 14 & 32 & 16 &  & \\
0.157 & \multicolumn{2}{c}{\ \dag} & $\infty$ & Triplicate \textbf{4.7} & Allemang \cite{allemang_2002} \\
0.2\.4\"1 & 10 & 4 & 36 & Cousins of \bfseries{0.317} & \\
0.31\"5 & 10 & 4 & 36 &  & \\
0.351 & 8 & 4 & 22 &  & Plambeck \\
0.51\.2 & 6 & 16 & 8 & $\Q \cong \mathcal{R}_8$ & \\
0.6\=4\=4 & 442 & 3255 & 172 &  & \\
0.71\.2 & 6 & 3 & 14 & $\Q \cong \mathcal{R}_{14}$ & \\
0.71\.6 & 2 & 22 & 14 & $\Q \cong \mathcal{R}_{14}$ & \\
4.5\.6 & 4 & 11 & 8 & $\Q \cong \mathcal{R}_8$ & \\
4.7\.4 & 2 & 8 & 8 & $\Q \cong \mathcal{R}_8$ & \\
\hline
\end{tabular}

\vspace{0.25cm}

$\dot d$ means ``$d$ or $d+1$''; $\ddot d$ means ``$d$ or $d+2$''; $\overline d$ means ``$d$, $d+1$, $d+2$ or $d+3$''

\vspace{0.25cm}
\dag\ Algebraic periodic; see Section~\ref{allemang} and the supplement~\cite{plambeck_2007supp}.
\caption{Wild two- and three-digit octal games with known mis\`ere quotients.}
\label{figure:solvedwildoctals}
\end{figure}

Comments and historical notes on some individual games follow.

\minisection{Guiles}
\label{guiles}
The game \textbf{0.15} is named \emph{Guiles} in honor of Richard Guy (the name is short for ``Guy's Kayles'').  Its mis\`ere quotient was discussed in Section~\ref{section:miserequotients}; it has 42 elements, with a $\SP$-portion of size 12.  See Figure~\ref{figure:015}, and also \cite{plambeck_200X}.

\minisection{Kayles}
\label{kayles}
William L. Sibert tells the interesting story of his discovery of the complete analysis of mis\`ere Kayles ({\bf 0.77})
in~\cite{sibert_1989}.  In~\cite{conway_1992}, Sibert's original solution is reformulated and simplified considerably.
The mis\`ere quotient of Kayles, a monoid with 40 elements, 
is discussed at length in \cite[Section 11.5]{plambeck_2005} and \cite{plambeck_200X}.
For completeness, we reproduce the solution in the supplement~\cite{plambeck_2007supp}.  Other octals with
mis\`ere quotient isomorphic to Kayles include {\bf 0.074}-{\bf 0.077} and {\bf 0.440}-{\bf 0.443} (duplicate Kayles), as well as {\bf 0.044} and {\bf 0.046} (triplicate Kayles).

\minisection{Allemang's Games}
\label{allemang} 
Allemang gives solutions for the wild octals {\bf 0.26}, {\bf 0.53}, {\bf 0.72}, {\bf 0.75} and {\bf 4.7}
in terms of his {\em generalized genus theory}~\cite{allemang_1984,allemang_2002}.  \textbf{0.53} and \textbf{0.72} have quotients of orders~16 and~24, respectively.  $\Q(\textbf{0.75}) \cong \mathcal{R}_8$; see Section~\ref{section:wildquotient}.  \textbf{0.26} and \textbf{4.7} are \emph{algebraic periodic} in the vague sense described in Section~\ref{section:algebraicperiodicity}; see Figures~\ref{figure:026} and~\ref{figure:47} for their presentations, and~\cite{plambeck_2007supp} for correctness proofs.

Allemang's solution to \textbf{0.26} is slightly flawed\footnote{For example, $H_{13} + H_{17} + H_{31}$ is a \PPos, but is misidentified in \cite{allemang_1984} as an \NPos.}; we give a corrected analysis in~\cite{plambeck_2007supp}.  His solution to \textbf{0.54} is also incorrect, but this appears somewhat more difficult to repair.  See Section~\ref{allemang-trouble} for further discussion of \textbf{0.54}.

\textbf{4.7} has several duplicates and triplicates; these are listed in Figure~\ref{figure:solvedwildoctals}.  In addition, the games \textbf{0.51\.6} and \textbf{0.57\.4} appear to have the same mis\`ere quotient as \textbf{4.7}, but we do not include solutions here.

\minisection{Flanigan's Games}
\label{flanigan}
\label{dot711}
Complete analyses of mis\`ere {\bf 0.34} and {\bf 0.71} are due to Jim Flanigan~\cite{berlekamp_1982}.
The game {\bf 0.34} has period eight (in both normal and mis\`ere play) and quotient order 12.
Although the normal-play Grundy sequence of {\bf 0.71} has period two, its mis\`ere period is six.  $|\Q(\textbf{0.71})| = 36$.

\minisection{Lemon Drops}
\label{lemondrops}
The game \textbf{0.56} is called \emph{Lemon Drops} in \emph{Winning Ways}~\cite{berlekamp_1982}.  It's tame, so the normal-play period of 144~\cite{gangolli_1989} remains the same in mis\`ere play.

\minisection{{\bf 0.123}}
\label{dot123}
This game is studied exhaustively in \cite{plambeck_2005}.  Its mis\`ere quotient has order~20.

\minisection{{\bf 0.241}}
\label{dot241}
The normal-play period of \textbf{0.241} is 2, while its mis\`ere-play period is~10.  \textbf{0.71}, discussed above, exhibits similar behavior.  Note that by Lemma~\ref{lemma:nkhconsequences}, the normal period of any faithfully regular octal game always divides its mis\`ere period.

\minisection{{\bf 0.644}}
\label{dot644}
This game has mis\`ere period 442 and preperiod 3255, the same as in normal play.  This is the largest known mis\`ere period for a wild octal game.
$\Q(\textbf{0.644})$ has order $172$.  The mis\`ere quotient ``grows'' for the last time at heap 333---that is, $\Q = \Q_{333}$, while $\Q \neq \Q_{332}$.  Intriguingly, the (normal) Grundy function also attains its maximum $\mathscr{G}$-value of 64 at heap 333.  This behavior is explained in the forthcoming paper~\cite{siegel_200Xe}.

\minisection{\textbf{0.4107}}
\begin{figure}
{\small
\pretend{12}{1 a a b b ab c c d e f g h b i ab^2 j 
k l m n o p q r abo anq b^3 s t abm cq^2 u cjk 
v w x b^3 y agt z b^2i \alpha{} \beta{} b^3 \gamma{} \delta{} b^4c bco ab\zeta{} \varepsilon{} 
b^3 \zeta{} ab^3c grx \eta{} ab\zeta{} ab\zeta{} b^3 ab^3c \theta{} b^3 b^4c ab^4 
cf\theta{} b^4c ab^3c ab^3c b^3 b^3 ab^4 b^4c b^4c ab^4 ab^3c 
b^3 b^3 ab^3c b^4c b^4c ab^4 ab^4 b^3 ab^3c ab^3c b^3 
b^4c ab^4 ab^4 b^4c ab^3c ab^3c b^3 b^3 ab^4 b^4c b^4c 
ab^4 ab^3c b^3 b^3 ab^3c b^4c b^4c ab^4 ab^4 b^3 ab^3c 
ab^3c b^3 b^4c ab^4 ab^4}
}
\caption{\label{4107sol} The values of {\bf 0.4107} have eventual period 24, despite many initial irregularities.}
\end{figure}
A great variety of mis\`ere quotients can be found among the four-digit octals; we include here just one particularly striking example.
\textbf{0.4107} has period~24, preperiod~66, and quotient order~506.  Its quotient has a minimal set of 34 generators  
\[\{a,b,c,\ldots,x,y,z,\alpha,\beta,\gamma,\delta,\varepsilon,\zeta,\eta,\theta\}.\]
The $\Phi$-values are shown in Figure \ref{4107sol}.  There are many irregular values among the smaller heaps, until finally a pattern abruptly emerges at heap size~66.  (See~\cite{miseregames_www} for the full quotient presentation.)

The sudden emergence of periodic behavior, after so much irregularity, is extraordinary.  One wonders how many other solutions lurk just beyond the reach of our computational resources.

\subsection{Unsolved Two- and Three-Digit Octals}
\suppressfloats[t]

We now briefly discuss some of the most important unsolved octal games.  Figure~\ref{figure:unsolvedtwodigits} lists the normal-play period of each unsolved two-digit octal.  Also listed is the largest $n$ for which the partial quotient $\Q_n$ is known, together with the order $|\Q_n|$.

\begin{figure}
\centering
\begin{tabular}{|>{\bfseries}l|rrrl|}
\hline
\mdseries Code & \multicolumn{1}{c}{N pd} & $n$ & $|\Q_n|$ & Comments \\
\hline
0.04 & \multicolumn{1}{c}{---} & 44 & 864 & Treblecross, \S \ref{treblecross} \\
0.06 & \multicolumn{1}{c}{---} & 15 & 48 & \\
0.\.07 & 34 & 33 & 638 & Dawson's Kayles, \S \ref{dawson} \\
0.14 & \multicolumn{1}{c}{---} & 20 & 96 & \\
0.16 & 149459 & 17 & 434 & \\
0.35 & 6 & 35 & 3182 & \\
0.36 & \multicolumn{1}{c}{---} & 20 & 304 & \S \ref{dot036} \\
0.37 & \multicolumn{1}{c}{---} & 15 & 304 & \S \ref{dot037} \\
0.4\=0 & 34 & 34 & 638 & Cousins of \bfseries{0.07} \\
0.4\"5 & 20 & 26 & 550 & \\
0.54 & 7 & \multicolumn{2}{c}{\dag} & \S \ref{allemang-trouble} \\
0.6\=0 & \multicolumn{1}{c}{---} & 16 & 304 & Officers \cite{berlekamp_1982} \\
0.6\=4 & \multicolumn{1}{c}{---} & 13 & 346 & \\
0.74 & \multicolumn{1}{c}{---} & 14 & 74 & \\
0.76 & \multicolumn{1}{c}{---} & 11 & 34 & \\
\hline
\end{tabular}

\vspace{0.25cm}

$\dot d$ means ``$d$ or $d+1$''; $\ddot d$ means ``$d$ or $d+2$''; $\overline d$ means ``$d$, $d+1$, $d+2$ or $d+3$''

\vspace{0.25cm}

\dag\ \textbf{0.54} appears to be algebraic periodic.
\caption{Unsolved two-digit octal games.}
\label{figure:unsolvedtwodigits}
\end{figure}

We suspect that every unsolved two-digit octal has an infinite mis\`ere quotient.  \textbf{0.54} is likely to be algebraic periodic in the sense of Section~\ref{section:algebraicperiodicity}; in the remaining cases, we suspect that the listed value of $n$ represents the last finite partial quotient.  Validating or refuting these suspicions will require more sophisticated techniques (or extremely diligent effort).

The unsolved three-digit octals are too numerous to list in a table of this form; they are summarized compactly in the supplement~\cite{plambeck_2007supp}, and on our website~\cite{miseregames_www} in more detail.  Some specific comments follow.

\minisection{Treblecross}
\label{treblecross}
{\bf 0.04} is a cousin of Treblecross, the game of ``one-dimensional Tic-Tac-Toe'' \cite{berlekamp_1982}.

\minisection{Dawson's Kayles}
\label{dawson}
Guy and Smith \cite{guy_1956} first observed that Dawson's Chess is equivalent to the octal game \textbf{0.137}.  It is a cousin of the two-digit octal \textbf{0.07}, which is commonly known as \emph{Dawson's Kayles}.

Guy and Smith showed that normal-play Dawson's Kayles (and therefore Dawson's Chess as well) has a period 34 Grundy
sequence.   Unpublished work by Ferguson~\cite{ferguson_unpub}, based on Conway's genus theory,
analyzed mis\`ere play of {\bf 0.07} to heap size 24.  Using {\em MisereSolver}, we can extend the analysis to heap size~33: the partial quotient $\Q_{33}(\mathbf{0.07})$ has order 638, with a $\SP$-portion of size 109.  The full presentation is too messy to justify its inclusion here, but it can be found online~\cite{miseregames_www} or reproduced with \emph{MisereSolver}.

The games \textbf{0.4\=0\=0} (in the schema notation of Figure~\ref{figure:unsolvedtwodigits}) are all equivalent to Dawson's Kayles.


\minisection{{\bf 0.36} and {\bf 0.37}}
\label{dot036}
\label{dot037}
The largest known partial quotients of {\bf 0.36} and {\bf 0.37} are isomorphic: $\Q_{20}({\bf 0.36}) \cong \Q_{15}({\bf 0.37})$.  However, their $\Phi$-values are somewhat different.  To what extent does this similarity continue?

\minisection{{\bf 0.54}}
\label{allemang-trouble}
Allemang \cite{allemang_2002} gives an incorrect solution to the game {\bf 0.54}.  It appears to be algebraic periodic, with an infinite mis\`ere quotient.  However, its quotient seems to be more complicated than those of {\bf 0.26} and {\bf 4.7}; better techniques for identifying algebraic periodicity are needed.  The three-digit octal game \textbf{0.145} exhibits similar (but not identical) behavior.


\minisection{{\bf 0.316}}
\label{dot316}
$|\Q_{23}(\textbf{0.316})| = 8704$.  This is the largest known finite partial quotient of an unsolved three-digit octal.

\minisection{{\bf 0.414}}
\label{dot414} 
The partial quotients $\Q_n(\textbf{0.4\"1\.4})$ grow surprisingly slowly as $n$ increases.  For this reason, \emph{MisereSolver} can quickly compute many corresponding $\Phi$-values.  Other wild games with similar behavior include \textbf{0.6\=4\=4}, \textbf{0.76\=4}, \textbf{0.77\.6}, and \textbf{4.4\=4}.  Among these, only \textbf{0.6\=4\=4} has been solved.  Since none of the others have known \emph{normal}-play solutions \cite{flammenkamp_www_octal}, it seems unlikely that their mis\`ere-play solutions are forthcoming.

Nonetheless, it is interesting to study this type of behavior.  There is a close relationship between slow-growing partial quotients and small $\SP$-portions; see~\cite{siegel_200Xe} for details.

\suppressfloats[t]
\subsection{Quaternary Games}

\begin{figure}
\centering
\begin{tabular}{@{}cc@{}}
\begin{tabular}{|c|ccc|}
\hline
Code & pd & ppd & $|\Q|$ \\
\hline
\rule{0pt}{11pt}{\bf 0.012\.2} & 7 & 8 & 20 \\
{\bf 0.1023} & 7 & 6 & 20 \\
{\bf 0.1032} & 7 & 8 & 20 \\
{\bf 0.1033} & 7 & 7 & 20 \\
{\bf 0.1231} & 5 & 5 & 20 \\
{\bf 0.123\.2} & 6 & 6 & 46 \\
{\bf 0.1321} & 5 & 6 & 20 \\
{\bf 0.1323} & 6 & 7 & 46 \\
{\bf 0.1331} & 5 & 5 & 20 \\
\hline
\end{tabular}
&
\begin{tabular}{|c|cccr|}
\hline
Code & pd & ppd & $|\Q|$ & Comments \\
\hline
\rule{0pt}{11pt}{\bf 0.2\.012} & 5 & 4 & 20 & \\
{\bf 0.3101} & 2 & 5 & 14 & \\
{\bf 0.3102} & \multicolumn{3}{c}{\!---} & $|\Q_{11}| =\ \,74$ \\
{\bf 0.3103} & 5 & 3 & 20 & \\
{\bf 0.3112} & 5 & 6 & 20 & \\
{\bf 0.3122} & \multicolumn{3}{c}{\!---} & $|\Q_6| =\ \,52$ \\
{\bf 0.3123} & \multicolumn{3}{c}{\!---} & $|\Q_{11}| = 328$ \\
{\bf 0.3131} & 2 & 7 & 12 & \\
{\bf 0.3312} & \multicolumn{3}{c}{\!---} & $|\Q_{13}| = 264$ \\
\hline
\end{tabular}
\end{tabular}
\caption{The twenty-one wild four-digit quaternary games.}
\label{figure:quaternary_summary}
\end{figure}

A \emph{quaternary game} is an octal game whose code digits are restricted to \textbf{0}, \textbf{1}, \textbf{2} and \textbf{3}.  Thus heaps may never be split, and the only available moves are to remove some number of tokens from a heap.  Despite these severe rules restrictions, many quaternary games exhibit interesting and surprisingly intricate mis\`ere quotients.

There are twenty-one wild four-digit quaternaries.  Of these, seventeen have known solutions; these are summarized in Figure \ref{figure:quaternary_summary}.  Only \textbf{0.3102}, \textbf{0.3122}, \textbf{0.3123}, and \textbf{0.3312} remain unsolved.  The first author has offered a reward of \$200 for the solution to \textbf{0.3102}, and \$25 each for the others.  It appears likely that their mis\`ere quotients are infinite, so the solutions might be quite difficult to obtain.  Figure~\ref{figure:quaternary_summary} lists their largest known partial quotients.

\suppressfloats[t]
\subsection{Algebraic Periodicity}
\label{section:algebraicperiodicity}

According to the Periodicity Theorem (Fact~\ref{fact:periodicitytheorem}), the full mis\`ere quotient $\Q(\Gamma)$ can sometimes be obtained from a finite number of partial quotients.  In order to apply the Perodicity Theorem, the partial quotients must ``stabilize'' after a certain point, so that $\Q_n(\Gamma) = \Q(\Gamma)$ for all sufficiently large~$n$.

Many mis\`ere games show another, more intriguing type of limiting behavior.  Such games have progressively larger partial quotients that exhibit a strong algebraic regularity.  When this regularity continues indefinitely, the full quotient can be deduced from a finite number of partial quotients; but unlike in the ``stable'' case, the full quotient is infinite, whereas every partial quotient is finite.  We've christened this behavior ``algebraic periodicity,'' but we can't give a precise definition because we don't fully understand how to describe it in~general.

\begin{figure}
\centering
\begin{tabular}{l@{}c}
\begin{minipage}{6.75cm}
$\begin{array}{c@{~}l}
\mathcal{Q} \cong \langle a,b,c_n~| &
{a^2=1},\ {b^{n+1}c_n=b^{2n+3}},\\
& {(c_mc_n=b^{m+2}c_n)_{m \leq n}} \rangle
\end{array}$

\vspace{0.5cm}

$\mathcal{P} = \{a,(b^{2m})_{m \geq 1},(b^mc_n)_{m \leq n\ \textrm{and}\ m+n\ \textrm{odd}}\}$
\end{minipage}
&
\pretend{4}{1 a b ab 1 a b ab c_0 ac_0 c_1 ab^3 c_2 abc_1 c_3 abc_2 c_4 abc_3 c_5 abc_4 c_6 abc_5}
\end{tabular}

\caption{\label{figure:026}Presentation and quotient map for \textbf{0.26}.}
\bigskip\medskip

\begin{tabular}{l@{}c}
\begin{minipage}{8.25cm}
\centering
$\begin{array}{c@{~}l}
\mathcal{Q} \cong \langle a,b,c,d_n~| &
{a^2=1},\ {bc=ab^3},\ {c^2=b^4}, \\
& {b^{n+1}d_n=a^{n+1}b^{2n+5}}, {cd_n=ab^2d_n}, \\
& {d_md_n=a^{m+1}b^{m+4}d_n} \rangle
\end{array}$

\vspace{0.5cm}

$\begin{array}{c@{}l}
\P = \{ & a,(b^{2m})_{m \geq 1},(b^md_n)_{m\ \textrm{odd},\ n\ \textrm{even},\ m < n},\\
& (ab^md_n)_{m\ \textrm{even},\ n\ \textrm{odd},\ m < n}\}
\end{array}$

\end{minipage}
&
\pretend{2}{a b a b c b^3 d_0 d_1 d_2 d_3 d_4 d_5 d_6}
\end{tabular}

\caption{\label{figure:47}Presentation and quotient map for \textbf{4.7}.}
\end{figure}

Two examples, the octal games \textbf{0.26} and \textbf{4.7}, are shown in Figures~\ref{figure:026} and~\ref{figure:47}.  Since we do not have any computational methods for verifying algebraic periodicity, we must resort to manual proofs of these figures.  The proofs are unenlightening, so we've relegated them to the supplement~\cite{plambeck_2007supp}.

Algebraic periodicity is a rich area for further study.  Many other games show this type of behavior, including \textbf{0.54}, \textbf{0.145}, \textbf{0.157}, \textbf{0.175}, \textbf{0.355}, \textbf{0.357}, \textbf{0.516}, \textbf{0.724}, and \textbf{0.734}.  We suspect that a general method for identifying and generalizing algebraic-periodic behavior---if one can be found---would quickly dispense solutions to many of them.

\begin{problem}
Give a precise definition of algebraic periodicity (that includes at least some of the cases listed above), and prove an analogue of the Periodicity Theorem for algebraic periodic games.
\end{problem}

\suppressfloats[t]
\subsection{Elements with Unusual Periods}

\label{subsection:unusualperiods}

If $\Q$ is a monoid and $x \in \Q$, then the \emph{period} of $x$ is the least $k \geq 1$ such that $x^{n+k} = x^n$, for some $n$.  If no such $k$ exists, then we say that $x$ has period $\infty$.  Note that if $\Q$ is finite, then $x$ must have finite period.

In normal play, every position $G$ satisfies the equation $G + G = 0$.  Therefore, in a normal quotient, every element has period at most two.  For a long time, we believed that the same is true for mis\`ere quotients.  However, there do exist finite mis\`ere quotients with larger-period elements, though they are exceedingly rare.

\begin{figure}
\centering
\includegraphics{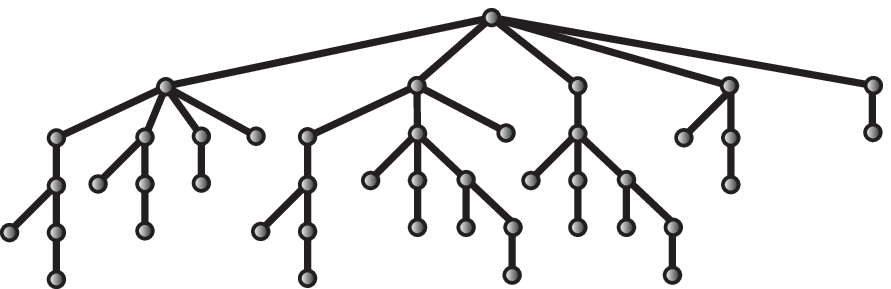}
\caption{\label{gametree}Ready for a game of mis\`ere coin-sliding on $G = \Star(2_\sh 210)(2_\sh 30)3_\sh 21$.}
\end{figure}

A striking example is given by $G = \Star(2_\sh 210)(2_\sh 30)3_\sh 21$.  It can be described as a coin-sliding game with heaps of tokens placed on the vertices of the tree shown in Figure~\ref{gametree}.  In this game, the players take turns sliding a single coin ``down'' the tree along a single edge.  The game ends when all coins have reached leaf nodes of~$G$, and whoever makes the last move loses.

$\Q(G)$ is a mis\`ere quotient of 200 elements; it has the following presentation:

\[
\Q(G) \cong \hangingpresent[t]{6cm}{a,b,c,d,e}{a2=1,b4=b2,b2c2=c2,c3=abc2,c2d=ac2,d2=c2,c2e2=c2,e6=1}
\]

\[
\begin{array}{c@{}l}
\P = \{ & a,b^2,abc,ab^3c,c^2,bd,b^3d,ae,ab^2de,cde,b^3cde,ae^2,b^2e^2,abce^2,ab^3ce^2, \\
& bde^2,ae^3,bde^3,b^3de^3,ae^4,b^2e^4,abce^4,ab^3ce^4,ab^2de^4,bcde^4,b^3cde^4,ae^5\}
\end{array}
\]

The most striking feature of this quotient is the generator $e$ of period 6.  In particular, $e^2$ has period 3, an odd number.

There also exist quotients with period-4 elements.  For example, let $H = \Star(2_\sh 1)(2_\sh 0)4310$.  Then $\Q(H)$ is a quotient of order 120 containing an element of period 4.  It is the smallest known quotient containing an element with period larger than 2.

\begin{problem}
Exhibit a mis\`ere quotient containing an element of period 8 (or any finite period other than 1, 2, 3, 4, or 6).
\end{problem}

\begin{question}
Are there any restrictions on the possible periods of mis\`ere quotient elements?
\end{question}

\begin{question}
What is the \emph{smallest} quotient containing an element of period $> 2$?
\end{question}

\subsection{The Sad Affair of $\Star(2_\sh 0)0$}
\label{subsection:infinitequotient}
\suppressfloats[t]

We now show that $\mathcal{Q}(\Star(2_\sh 0)0)$ is infinite.  It can be verified computationally that all games with smaller game trees have finite quotients, so in some sense $\Star(2_\sh 0)0$ is the \emph{simplest} game that yields an infinite quotient.  This fact makes it worthy of significant attention.

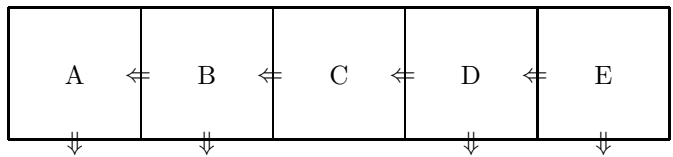
\begin{figure}
\centering
\begin{picture}(250,50)
\put(0,0){\line(1,0){250}}
\put(0,50){\line(1,0){250}}
\multiput(0,0)(50,0){6}{\line(0,1){50}}
\multiput(44,22)(50,0){4}{\makebox{$\Leftarrow$}}
\multiput(22,-4)(50,0){2}{\makebox{$\Downarrow$}}
\multiput(172,-4)(50,0){2}{\makebox{$\Downarrow$}}
\put(20,21){\makebox[10pt][c]{A}}
\put(70,21){\makebox[10pt][c]{B}}
\put(120,21){\makebox[10pt][c]{C}}
\put(170,21){\makebox[10pt][c]{D}}
\put(220,21){\makebox[10pt][c]{E}}
\end{picture}
\caption{A simple coin-sliding game with an infinite mis\`ere quotient.}
\label{figure:coinsliding_inf_quotient}
\end{figure}

$\Star(2_\sh 0)0$ arises in various quaternary games; for example, it is the canonical form of an $\mathbf{0.31011}$-heap of size~5.  It can also be represented as the coin-sliding game shown in Figure~\ref{figure:coinsliding_inf_quotient}.  Stacks of coins are arranged on the boxes to form a starting position.  On her turn, a player may slide one coin one space in the direction of an arrow.  This might cause the coin to drop off the board, whereupon it is removed from the game.  It is easily seen that a single coin on the rightmost box has canonical form $\Star(2_\sh 0)0$.

For the remainder of this section, write
\[A = \Star;\quad B = \Star 2;\quad C = \Star2_\sh ;\quad D = \Star2_\sh 0;\quad E = \Star(2_\sh 0)0;\]
and put $\mathscr{A} = \mbox{cl}(E)$.  Then each game in $\mathscr{A}$ has the form
\[iA + jB + kC + lD + mE,\]
for some quintuple of integers $(i,j,k,l,m)$.  (These correspond to the number of coins on each successive box in Figure \ref{figure:coinsliding_inf_quotient}, left-to-right.)

When $k \geq 3$, the \PPos{}s admit a simple description: $G$ is a \PPos{} iff both $i+l$ and $j+m$ are even.  When $k \leq 2$, however, the \PPos{}s are highly erratic.  The outcomes for $l < 14$ and $m < 18$ are summarized in Figure~\ref{figure:ppos_inf_quotient}.  Each of the twelve grids represents the outcomes for a particular triple $(i,j,k)$.  Within each grid, there is a black dot at (row $m$, column $l$) iff $G = iA + jB + kC + lD + mE$ is a \PPos.  The outcomes for $l \geq 14$ and $m \geq 18$ can be obtained through the following simple recurrence:

\begin{figure}
\centering
\begin{tabular}{@{}c@{\hspace{-8pt}}c@{\hspace{-8pt}}c@{\hspace{-8pt}}c@{}}
& \hspace{8pt} $k=0$ & \hspace{8pt} $k=1$ & \hspace{8pt} $k=2$ \vspace{4pt} \\
\input ppostable.tex
\end{tabular}

\[\bullet = \SP; \qquad \textrm{x} = \SP \textrm{ iff } j \geq 2; \qquad \circ = \SP \textrm{ iff } j < 2\]

\vspace{-11pt}

\caption{Schematic of the \PPos{}s for $\cl(\Star(2_\sh 0)0)$ with $k \leq 2$.}
\label{figure:ppos_inf_quotient}
\end{figure}

\begin{center}
Let $G = iA + jB + kC + lD + mE$, with $l \geq 4$ or $m \geq 12$.  Then:
\end{center}
\[
\begin{array}{r@{~=~}ll}
o^-(G + 2D + 2E) & o^-(G); \\
o^-(G + 2D) & o^-(G) & \textrm{if $m < 3$}; \\
o^-(G + 2E) & o^-(G) & \textrm{if $l < 3$}.
\end{array}
\tag{\dag}
\label{eq:recurrence_inf_quotient}
\]

Figure~\ref{figure:ppos_inf_quotient} illustrates some striking features about this game.  Most prominent is the diagonal of \PPos{}s in the $(0,0,0)$ diagram, along the line $m=l+7$.  In fact, it is this diagonal (and its echo in several other cases) that makes the quotient infinite.  In addition, there are many strange anomalies, such as the \PPos{} at $(i,j,k,l,m)=(1,1,2,0,3)$.  Such anomalies are difficult to explain as anything other than combinatorial chaos.  It is remarkable that such a simple game gives rise to so much complexity.

A comparison with normal play is instructive.  It is easily checked that $A,B,C,D,E$ have Grundy values $1,2,0,1,2$, respectively.  Therefore~$G$ is a \emph{normal} \PPos{} iff both $i+l$ and $j+m$ are even.  Note that this coincides with the mis\`ere condition when $k \geq 3$.  This is explained by our discussion of normality in Section~\ref{section:nkh}: put $z = \Phi(3C)$; then $x \mapsto zx$ is a surjective homomorphism from $\Q(E)$ onto the group $\mathbb{Z}_2 \times \mathbb{Z}_2$.  Since~$z$ is an idempotent and~$z\Q$ is a group, $z$ must be the kernel identity of~$\Q$.  Thus we have a beautiful illustration of the generalized mis\`ere strategy described in Section~\ref{section:nkh}: follow normal play unless your move would leave fewer than three coins on box~$C$.





We now prove the correctness of Figure~\ref{figure:ppos_inf_quotient} (and of the asserted outcomes when $k,l,m$ are large).  For each $G \in \A$, denote by $o^*(G)$ the \emph{asserted} outcome of $G$; as always, $o^-(G)$ denotes the \emph{true} outcome of $G$.

\begin{lemma}
\label{lemma:ostar_recurrences}
Let $G = iA + jB + kC + lD + mE$.
\begin{enumerate}
\item[(a)] If $j \geq 2$, then $o^*(G+2B) = o^*(G)$.
\item[(b)] If $k \geq 3$, then $o^*(G+C) = o^*(G)$.
\end{enumerate}
\end{lemma}

\begin{proof}
(a) is by inspection of Figure~\ref{figure:ppos_inf_quotient}; (b) is automatic from the description of~$o^*$.
\end{proof}

\begin{theorem}
$o^-(G) = o^*(G)$ for all $G \in \A$.
\end{theorem}

\begin{proof}
Note that $o^*(0) = \mathscr{N}$.  Thus it suffices to show the following, for each $G \neq 0$:
\[o^*(G) = \mathscr{P} \Longleftrightarrow \textrm{No option } G' \textrm{ satisfies } o^*(G') = \mathscr{P}.
\tag{\ddag}\label{eq:ostar_sufficient}\]

Suppose instead (for contradiction) that (\ref{eq:ostar_sufficient}) fails for some $G \in \A$.  Choose $G = iA+jB+kC+lD+mE$ to be a counterexample with minimal birthday.  Since $A+A = 0$ (canonically), we know that $i = 0$ or~$1$.

Now (\ddag) can be verified computationally when $j < 5$, $k < 5$, $l < 7$, and $m < 15$.  Therefore we may safely assume that at least one of these inequalities fails.  There are five cases.

\vspace{0.15cm}\noindent
\emph{Case 1}: $j \geq 5$.  Then fix $H$ with $G = H + 2B$.  Since $j - 2 \geq 2$, Lemma~\ref{lemma:ostar_recurrences}(a) gives $o^*(G) = o^*(H)$.  Furthermore, every option $G'$ of $G$ can be written as $G' = H' + 2B$.  Since $j - 2 \geq 3$, we can write $H' = j'B + X$, for $j' \geq 2$.  Thus Lemma~\ref{lemma:ostar_recurrences}(a) gives $o^*(H' + 2B) = o^*(H')$.  Therefore
\[\{o^*(G') : G' \in \opts(G)\} = \{o^*(H') : H' \in \opts(H)\}.\]
It follows that (\ddag) fails for $H$ as well, contradicting minimality of $G$.

\vspace{0.15cm}
The remaining cases are all very similar.

\vspace{0.15cm}\noindent
\emph{Case 2}: $k \geq 5$.  Then fix $H$ with $G = H + C$.  Since $k - 1 \geq 4$, Lemma~\ref{lemma:ostar_recurrences}(b) gives $o^*(G) = o^*(H)$, but also $\{o^*(G') : G' \in \opts(G)\} = \{o^*(H') : H' \in \opts(H)\}$, contradicting minimality of $G$.

\vspace{0.15cm}\noindent
\emph{Case 3}: $l \geq 7$ but $m < 3$.  Then fix $H$ with $G = H + 2D$.  Since $l - 2 \geq 5$, the recurrence (\dag) gives $o^*(G) = o^*(H)$, etc., as in previous cases.

\vspace{0.15cm}\noindent
\emph{Case 4}: $m \geq 15$ but $l < 3$.  Then fix $H$ with $G = H + 2E$, and proceed as in Case~3.

\vspace{0.15cm}\noindent
\emph{Case 5}: $l \geq 7$ or $m \geq 15$, and $l,m \geq 3$.  Then fix $H$ with $G = H + 2D + 2E$, and proceed as in Case~3.
\end{proof}

\begin{corollary}
$\Q(E)$ is infinite.
\end{corollary}

\begin{proof}
If $l \geq 3$ is odd, then $lD + mE$ is a \PPos{} iff $m = l + 7$.  Thus the games $(2n+3) \cdot D$ are pairwise distinguishable.
\end{proof}

\subsection{Irregular and Abnormal Quotients}

In Section~\ref{section:nkh} we defined a quotient to be \emph{regular} if $|\mathcal{K} \cap \P| = 1$, and \emph{normal} if $\mathcal{K} \cap \P = \{z\}$.  Irregular quotients are difficult to find; the smallest known example is $\Q_{12}(\textbf{0.324})$, with $|\mathcal{K} \cap \P| = 4$.  Simpler examples almost certainly exist, but better techniques are needed in order to detect them.



\section*{Acknowledgements}

The authors wish to thank Dan Hoey for many helpful comments and suggestions, and the referees for their extremely careful and insightful reading of the manuscript.

\nocite{allemang_2002unpub}
\nocite{sprague_1947}
\nocite{siegel_mqlectures}
\bibliography{games}

\end{document}